\documentclass[a4paper,reqno,11pt]{amsart}
\usepackage{amsmath, latexsym, amsfonts, amssymb, amsthm, amscd}
\usepackage{graphics,epsf,epsfig,psfrag}

\usepackage[pdftex]{hyperref}

\numberwithin{equation}{section}

\newtheorem{theorem}{Theorem}[section]
\newtheorem{lemma}[theorem]{Lemma}
\newtheorem{proposition}[theorem]{Proposition}

\newtheorem{rem}[theorem]{Remark}

\newtheorem{assumption}[theorem]{Assumption}
\newtheorem{hypothesis}[theorem]{Hypothesis}

\newcommand{\bbE}{{\ensuremath{\mathbb E}} }

\newcommand{\bbP}{{\ensuremath{\mathbb P}} }

\newcommand{\bbT}{{\ensuremath{\mathbb T}} }

\newcommand{\cB}{{\ensuremath{\mathcal B}} }

\newcommand{\cF}{{\ensuremath{\mathcal F}} }

\newcommand{\cL}{{\ensuremath{\mathcal L}} }

\newcommand{\ga}{\alpha}
\newcommand{\gb}{\beta}

\newcommand{\gd}{\delta}

\newcommand{\gep}{\varepsilon}       

\newcommand{\gl}{\lambda}

\newcommand{\gs}{\sigma}

\newcommand{\gp}{\varphi}

\newcommand{\go}{\omega}

\renewcommand{\tilde}{\widetilde}          
\DeclareMathSymbol{\leqslant}{\mathalpha}{AMSa}{"36} 
\DeclareMathSymbol{\geqslant}{\mathalpha}{AMSa}{"3E} 
\DeclareMathSymbol{\eset}{\mathalpha}{AMSb}{"3F}     
\newcommand{\dd}{\text{\rm d}}             

\DeclareMathOperator*{\union}{\bigcup}       

\newcommand{\R}{\mathbb{R}}
\newcommand{\C}{\mathbb{C}}
\newcommand{\Z}{\mathbb{Z}}
\newcommand{\N}{\mathbb{N}}

\DeclareMathOperator{\var}{var}
\DeclareMathOperator{\cov}{cov}

\def\bP{\ensuremath{\bs{\mathrm{P}}}}
\def\bE{\ensuremath{\bs{\mathrm{E}}}}

\def\bs{\boldsymbol}
\def\proof{\noindent\textbf{Proof.} \ }

\renewcommand{\phi}{\varphi}

\DeclareMathOperator{\tr}{Tr}

\newcommand{\rbar}{\overline{r}}

\newcommand{\somb}{\widehat} 

\begin{document}

\title[Disordered semiflexible polymers]{Large scale behavior of semiflexible heteropolymers}

\author{Francesco Caravenna}

\address{Dipartimento di Matematica Pura e Applicata, Universit\`a
degli Studi di Padova,
via \mbox{Trieste} 63, 35121 Padova, Italy}

\email{francesco.caravenna@math.unipd.it}

\author{Giambattista Giacomin}
\address{
Universit\'e Paris  Diderot (Paris 7) and Laboratoire de Probabilit{\'e}s et Mod\`eles Al\'eatoires (CNRS U.M.R. 7599),
                             UFR de Math\'ematiques,
                             case 7012 (site Chevaleret)
                             75205 PARIS Cedex 13, France
}

\email{giacomin@math.jussieu.fr}

\author{Massimiliano Gubinelli}
\address{Equipe de probabilit\'es, statistique et mod\'elisation, Universit\'e de Paris-Sud, B\^atiment 425, 91405 ORSAY Cedex, France
-- now at CEREMADE, Universit\'e Paris Dauphine, France}

\email{gubinelli@ceremade.dauphine.fr}



\begin{abstract}
We consider a general discrete model for heterogeneous semiflexible polymer chains.
Both the thermal noise and the inhomogeneous
character of the chain (the {\sl disorder}) are modeled in terms of random rotations.
We focus on the {\sl quenched} regime, {\sl i.e.}, the analysis is performed 
for a given realization of the disorder. Semiflexible models differ substantially from 
random walks on short scales, but on large scales a Brownian behavior emerges.
By exploiting techniques from tensor analysis and non-commutative Fourier analysis,
we establish the Brownian character of the model on large scales
and we obtain an expression for the diffusion constant. We moreover give conditions
yielding quantitative mixing properties.

\bigskip
\medskip

\noindent
{\sc R\'esum\'e.}
On consid\`ere un mod\`ele discret pour un polym\`ere semi-flexible et h\'et\'erog\`ene. 
Le bruit thermique et le caract\`ere h\'et\'erog\`ene du polym\`ere (le {\sl d\'esordre}) 
sont mod\'elis\'es en termes de rotations al\'eatoires. Nous nous concentrons sur 
le r\'egime de d\'esordre {\sl g\'el\'e}, c'est-\`a-dire, l'analyse est effectu\'ee pour une 
r\'ealisation fix\'ee du d\'esordre. Les mod\`eles semi-flexibles diff\`erent sensiblement 
des marches al\'eatoires \`a petite \'echelle, mais \`a grande \'echelle un comportement 
brownien appara\^it. En exploitant des techniques de calcul tensoriel et d'analyse 
de Fourier non-commutative, nous \'etablissons le caract\`ere brownien du mod\`ele 
\`a grande \'echelle et nous obtenons une expression pour la constante de diffusion. 
Nous donnons aussi des conditions qui entra\^inent des propri\'et\'es quantitatives de m\'elange.
\end{abstract}

\date{\today}

\maketitle

\section{Introduction}

\subsection{Homogeneous semiflexible polymer models}
In the vast polymer modeling literature an important role is played
by random walks, in fact self-avoiding random walks ({\sl e.g.} \cite{cf:Freed,cf:scaling}).
However they are expected to model properly real polymers
only on large scales. On shorter scales one observes a stiffer behavior
of the chain,  and other models
have been proposed, notably  the {\sl semiflexible} one (see {\sl e.g.} \cite{cf:MarkoSiggia,cf:WN}
and references therein).
A {\sl semiflexible polymer} is a natural and appealing mathematical
object and, in absence of self-avoidance,
it has been implicitly considered in the probability  literature
for a long time. 
Consider in fact a probability measure $Q$ on the Lie group $SO(d)$ --
the rotations in $\R^d$ ($d=2,3,\ldots$) -- and sample from this, in an independent
fashion, a sequence of rotations $r_1, r_2, \ldots\,$.
Fixing an arbitrary rotation $R \in SO(d)$ and denoting by
$e^1, \ldots, e^d$ the unit coordinate vectors in $\R^d$,
the process $\{v_n\}_{n \ge 0}$ defined by
\begin{equation} \label{eq:rw}
v_0 \;:=\; R\, e^d \,, \qquad \quad
v_n \;:=\; \big(R\, r_1\, r_2\, \cdots \, r_n \big) e^d\,, \quad n=1, 2 , \ldots \,,
\end{equation}
is nothing but a {\sl random walk} on the unit sphere $S^{d-1} \subset \R^d$
starting at $v_0$, a much studied object ({\sl e.g.} \cite{cf:RobUrs,cf:Ros}).
Then the process $\{X_n^{v_0}\}_{n=0, 1, \ldots}$ defined by
\begin{equation}
\label{eq:Xhomogeneous}
X_n^{v_0} \; := v_0 \; + \sum_{j=1}^n v_j\, =\, \sum_{j=0}^n v_j,
\end{equation}
is a {\sl homogeneous} semiflexible polymer model in dimension~$d$.
The reason for writing $(R \, r_1\, r_2\, \cdots \, r_n)$ instead
of $(r_n\, r_{n-1}\, \cdots \, r_1 \, R)$ in \eqref{eq:rw} is explained
in Remark~\ref{rem:order} below.

\begin{figure}[h]
\begin{center}
\leavevmode
\includegraphics[width=12.8cm]{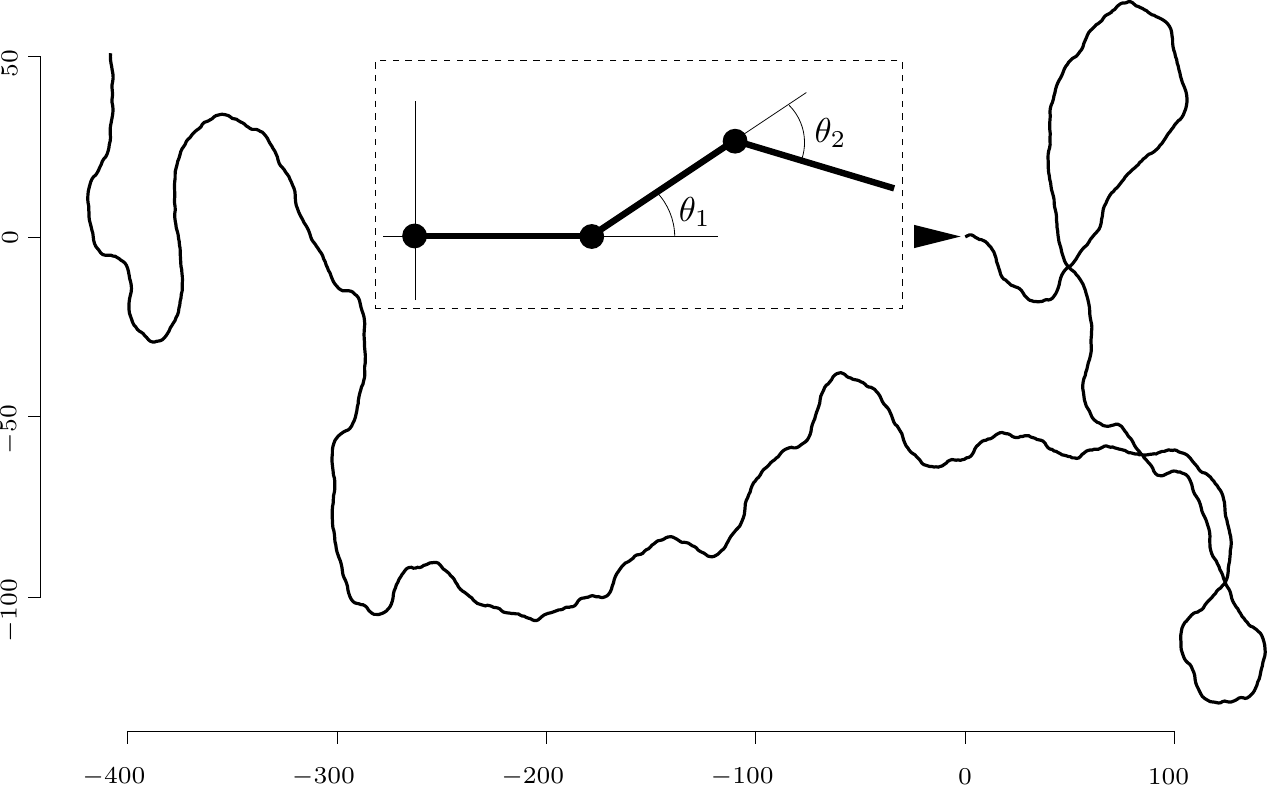}
\end{center}
\caption{\label{fig:1} \small 
A sample bidimensional trajectory, with $n=10^4$ and $\{\theta_i\}_{i=1,2, \ldots}$
drawn  uniformly from  $(-\pi/10, \pi/10)$, while $\theta_0$ is $0$ (the notation
is the one of Remark~\ref{rem:2dsimple}). In the inset there is a zoom of the
starting portion of the polymer (the starting point is marked by the arrow).
It is clear that the starting orientation 
$v_0=(1,0)$ sets up a drift that is forgotten only after a certain number 
of steps. Moreover, even if the starting orientation 
eventually fades away, in the sense that the expectation 
of the scalar product of $v_n$ and $v_0$ vanishes as $n$ becomes large, the
local orientation is carried along for a while. 
A precise meaning to this  is brought by the  
key concept of {\sl persistence length} $\ell$, that can be defined 
as the reciprocal of the rate of exponential decay of $\bE \langle v_0 , v_n \rangle$,
where $\langle \cdot , \cdot \rangle$ denotes the standard scalar product
in $\R^d$ and $\bE$ is the average over the variables $\{r_i\}_i$.
Intuitively,  one expects that on a scale much larger than the persistence length,
the semiflexible polymer $X_n^{v_0}$ is going to behave {\sl like
a random walk}. Note that if we view the 
elements of $SO(d)$ as linear operators, we can define 
$\rbar := \bE r_1$ (not a rotation unless $r_1$ is trivial!) and we have
$\bE \langle v_0 , v_n \rangle = \langle e^d , \rbar ^n e^d \rangle$,
which shows that the decay of $\bE \langle v_0 , v_n \rangle$ is indeed
of exponential type.
}
\end{figure}

\medskip
\begin{rem}\label{rem:2dsimple}
\rm The reader can get some intuition on the process by having a look at 
the two-dimensional case of Figure~\ref{fig:1}. This case
is in reality particularly 
easy to analyze in detail (and it does not capture the full
complexity of the $d>2$ case) because the rotations in two
dimensions commute and they 
are characterized by only one parameter.
More precisely, if we identify the random rotation $r_j$
with the angle $\theta_j$, for $j \ge 1$, and we take $\theta_0$
such that $v_0 = (\cos \theta_0, \sin \theta_0)$, by setting
$\gp_n:= \theta_0 + \theta_1 + \ldots +
\theta_n$ we can write
\begin{equation}
	X_n^{v_0} \; =\; v_0 \,+\, \left( \sum_{j=1}^n \cos(\gp_j) \,,\,  
	\sum_{j=1}^n \sin(\gp_j)\right) \,.
\end{equation}
This explicit expression allows an easy and complete
analysis of the two-dimensional case, {\sl cf.} Appendix~\ref{sec:2d}.
Of course, in general no such simplification is possible for $d > 2$.
\end{rem}

\smallskip

Homogeneous semiflexible chains  have been used in a variety
of contexts \cite{cf:WN,cf:ZC} and they do propose challenging questions that 
are still only partially understood (even in their continuum version,
see Remark~\ref{rem:continuum}),
also because it is difficult to obtain explicit expressions for very basic quantities
like the {\sl loop formation probability}, {\sl i.e.} the hitting probability. 
As a matter of fact, a more realistic model would have to take into account 
a self-avoiding constraint, which is more properly called 
{\sl excluded volume condition}, that imposes
that the {\sl sausage-like} trajectory does not self-intersect. This of course 
makes the model extremely difficult to deal with. 
Added to that,  models need to 
embody the fact that often real polymers are
{\sl inhomogeneous}, {\sl i.e.} they are not made up
of identical monomers and that this does affect the geometry of 
the configurations. It is  precisely on this latter direction that we are going 
to focus. 


\smallskip

\subsection{Heterogeneous models}
Heterogeneous semiflexible chains have attracted a substantial
amount of attention (see {\sl e.g.} \cite{cf:BDM,cf:MN,cf:Vaillant,cf:WN,cf:ZC}),  
often (but not only) as  a  modeling frame  for DNA or RNA (single or double stranded) chains. 
The information that we want to incorporate in the model
is the fact that the monomer units may vary along the chain:
for the DNA case, the four bases A, T, G and C are the origin
of the inhomogeneity and couple of monomer units  have an associated typical bend
that depends on their bases. The model we are interested
in is therefore still based on randomly sampled rotations $r_1, r_2, \ldots$,
independent and identically distributed with a given marginal law $Q$
(this represents the thermal noise in the chain), but associated to that there is a sequence
of rotations $\go_1, \go_2,\ldots$ that is fixed and does not
fluctuate with the chain. If we want to stick to the  DNA example,
the $\go$-sequence is fixed once the base sequence is given.
The model is then defined by giving once again the orientation
$v_0 = R \, e^d \in S^{d-1}$
of the initial monomer and by defining for $n \ge 0$
\begin{equation} \label{eq:X}
X_n^{v_0, \go} \; := \; v_0 \,+\, \sum_{j=1}^n v_j^\go \qquad \text{ with } \qquad
v_j^\go\, :=\,  \big( R \,\go_1\, r_1 \cdots \go_j \, r_j \big) e^d \,.
\end{equation}
It should be clear that the rotation $\go_i$ sets up the equilibrium position
of the $i$-th monomer with respect to the $(i-1)$-st. 
In different terms, the
sequence $\go_1, \go_2, \ldots$ defines the {\sl backbone} around which the semiflexible
chain fluctuates.
 
\smallskip
 
The aim of this paper is to study 
the large scale behavior of the process $\{X_n^{v_0, \go}\}_{n}$ when the sequence
$\go$ is {\sl disordered}, {\sl i.e.} it is chosen as the typical
realization of a random process. The simplest example
is of course the one in which the variables $\go_n$ are independent
and identically distributed, but we stress from now that
we are interested in the much more general case when $\go$ is an
{\sl ergodic process} (see Assumption~\ref{ass:disorder} for
the definition of ergodicity). This includes strongly correlated 
sequences of random variables and, in particular, the ones that
have been proposed to mimic the base distributions along the DNA
({\sl e.g.} \cite{cf:Vaillant} and references therein).
Other aspects of this model deserve attention, notably the analysis 
of the {\sl persistence length} in the heterogeneous set-up
(see the caption of Figure~\ref{fig:1})
and other kind of scaling limits, like the {\sl Kratky-Porod} limit
(see Remark~\ref{rem:continuum}):
these issues are taken up in a companion paper.


\medskip

\begin{rem}\label{rem:order} \rm
Let us comment on the order of the rotations appearing in equations \eqref{eq:X}
and \eqref{eq:rw}. The key point is the following consideration:
in defining the rotations $r_i$ and $\omega_i$, we assume that the
$i$-th monomer lies along the direction $e^d$. Therefore, before applying
these rotations, we have to express them in the actual reference frame of
the $i$-th monomer. 
Let us be more precise, considering first
the homogeneous case given by \eqref{eq:rw}.
The rotation $r_1$ describes the thermal
fluctuations of the first monomer {\sl assuming that its equilibrium position is $e^d$}.
However the equilibrium position of the first monomer is rather $v_0 = R\, e^d$, therefore
we first have to express $r_1$ in the reference frame of $v_0$,
obtaining $R \, r_1 \, R^{-1}$, and then apply it to $v_0$, obtaining
$v_1 = (R \, r_1 \, R^{-1})\, v_0 = (R \, r_1) \, e^d$. The same procedure
yields $v_2 = (R \, r_1) \, r_2 \, (R \, r_1)^{-1} \, v_1 = (R \, r_1 \, r_2) \, e^d$,
and so on. The inhomogeneous case of equation \eqref{eq:X} is analogous:
we first apply $\go_1$ expressed in the reference frame of $v_0$, getting
$v_0' = (R \, \go_1 \, R^{-1}) \, v_0 = (R \, \go_1) \, e^d$, then we apply
$r_1$ expressed in the reference frame of $v_0'$, obtaining
$v_1^\go = (R \, \go_1) \, r_1 (R \, \go_1)^{-1} \, v_0' = (R \, \go_1 \, r_1) \, e^d$,
and so on.
\end{rem}

\smallskip

\begin{rem}\label{rem:continuum}
\rm
Most of the physical literature focuses 
on 
a continuum version of the homogeneous 
semiflexible model, often called {\sl wormlike chain}  or {\sl Kratky-Porod} model
({\sl e.g.} \cite{cf:WN} and references therein),
which can be  obtained in a large scale/high stiffness limit
of discrete models.
As for the discrete semiflexible model we had a discrete length parameter
$n$ that was in fact counting the monomers along the chain, here we have 
a continuous parameter $t\ge 0$ and the location $\tilde X_t$ of 
the wormlike chain at $t$ is equal to 
$\int_{0}^t {\tt B}^{(d)}(s) \dd s$,
where $\left\{ {\tt B}^{(d)}(s)\right\}_{s\ge 0}$
is a Brownian motion on $S^{d-1}$ ({\sl e.g.} \cite{cf:RW}). Note that 
the initial orientation $v_0$ is here replaced by the choice of ${\tt B}^{(d)}(0)$.
For $d=2$, once again, this process becomes particularly 
easy to describe since ${\tt B}^{(2)}(t)= 
(\cos (B(t)+x_0), \sin(B(t)+y_0) )$,
where $B$ is a standard Brownian motion.
We point out that in the physical literature the continuum model is just used
for some {\sl formal} computations and, in the heterogeneous set-up, the model
is often ill-defined and in fact when simulations are performed usually one
goes back to a discrete model \cite{cf:BDM,cf:MN,cf:Vaillant,cf:WN,cf:ZC}.   
\end{rem}

\smallskip


\smallskip

\subsection{The Brownian scaling}

In order to study the large scale behavior of our model,
we introduce its diffusive rescaling, {\sl i.e.} the
continuous time process $B_{N}^{v_0, \go}(t)$ defined for $N \in \N$ 
and $tN\in \N\cup \{0\}$ by
\begin{equation} \label{eq:BN}
	B_{N}^{v_0, \go}(t) \; :=\;
	\frac{1}{\sqrt N} \, X^{v_0, \omega}_{ Nt } .
\end{equation}
This definition is extended to every 
$t \in [0, \infty)$ by linear interpolation,
so that $B_{N}^{v_0, \go}(\cdot)\in C([0,\infty))$ and it is piecewise affine,
where $C([0,\infty))$ denotes the space of
real-valued continuous functions defined on $[0,\infty)$ and is
equipped as usual with the topology of
uniform converge over the compact sets and with the
corresponding $\gs$-field. The precise hypothesis we make on
the thermal noise is as follows.

\medskip
\begin{assumption} \rm \label{ass:noise}
The variables $\big(\{r_n\}_{n\ge 1}, \bP \big)$ taking values in $SO(d)$ are independent
and identically distributed, and the 
law $Q$ of $r_1$ satisfies the following irreducibility condition:
there do not exist linear subspaces $V, W \subseteq \R^d$ such that
$Q(g \in SO(d): \, g\, V = W) = 1$, except the trivial cases when
$V = W = \{0\}$ or $V = W = \R^d$.
\end{assumption}

\smallskip

\noindent
We point out that this assumption on $Q$ (actually on its support) is very mild.
It is fulfilled for instance whenever the
support of $Q$ contains a non-empty open set $A \subseteq SO(d)$
(this is a direct consequence of the fact that an open subset of $SO(d)$
spans $SO(d)$), in particular when
$Q$ is absolutely continuous with respect to the
Haar measure on $SO(d)$, a very reasonable
assumption for thermal fluctuations (see \S\ref{sec:notations}
for details on the Haar measure). We stress however that 
absolute continuity is not necessary and in fact
several interesting cases of discrete laws are allowed (e.g.,
for $d=3$, when $Q$ is supported on the symmetry group of a Platonic solid).
Also notice that for $d=2$ Assumption~\ref{ass:noise} can be restated
more explicitly as follows: denoting by $R_{\theta} \in SO(2)$ the rotation by an
angle $\theta$, there does not exist $\theta \in [0,\pi)$ such
that $Q(\{R_{\theta}, R_{\pi + \theta}\})=1$.

\smallskip

Next we state precisely our assumption on the disorder.  

\medskip
\begin{assumption} \rm \label{ass:disorder}
The sequence $\big(\{\go_n\}_{n\ge 1},\bbP\big)$ is {\sl stationary},
{\sl i.e.} $\{\go_{n+1}\}_{n\ge 1} $ and $ \{\go_n\}_{n\ge 1}$ have the same 
law, and {\sl ergodic},
{\sl i.e.} $\bbP \big( \{\go_n\}_{n\ge 1} \in A \big) \in \{0,1\}$
for every shift-invariant measurable set $A \subseteq SO(d)^{\N}$.
Shift-invariant means that $\{x_1, x_2, \ldots\} \in A$
if and only if $\{x_2, x_3, \ldots\} \in A$, while measurability is with respect to the
product $\gs$-field on $SO(d)^{\N}$.
\end{assumption}
\medskip

We can now state our main result.

\medskip
\begin{theorem} \label{th:main}
If Assumptions~\ref{ass:noise} and~\ref{ass:disorder} are satisfied,
then $\bbP(\dd \go)$-almost surely and for every choice of $v_0$ the
process $B_{N}^{v_0, \go}$ converges in distribution on $C([0,\infty))$
as $N \to\infty$ toward $\gs B$, where 
$B=\left\{ \left(B_1(t), \ldots , B_d(t)\right) \right\}_{t\ge 0}$
is a standard $d$-dimensional Brownian motion and the positive
constant $\sigma^2$ is given by
\begin{equation}
\label{eq:sigma-def}
	\sigma^2 \;:=\; \frac 1 d \,+\, \frac 2 d\, \sum_{k=1}^\infty \,
	\bbE  \, \bE \, \langle e^d ,
	 \go_1 \, r_1 \, \cdots \, \go_k \, r_k  \, e^d \rangle \,,
\end{equation}
where the series in the right-hand side  converges.
\end{theorem}

\medskip

This result says, in particular, that the disorder affects the large 
scale behavior of the polymer only through the diffusion coefficient $\gs^2$.
Let us now consider some special cases in which $\gs^2$ can be made more explicit.
Notice first that, by setting $\rbar := \bE(r_1)$, we can rewrite
$\bbE  \, \bE \, \langle e^d , \go_1 \, r_1 \, \cdots \, \go_k \, r_k  \, e^d \rangle
= \bbE  \, \langle e^d , \go_1 \, \rbar \, \cdots \, \go_k \, \rbar \, e^d \rangle$.
\begin{itemize}

\item When $\overline r = c I$, where $I$
denotes the identity matrix and $c$ is a constant
(necessarily $|c| < 1$), the expression for $\gs^2$ becomes
\begin{equation}
	\gs^2 \;=\; \frac 1d \,+\, \frac 2d \,
	\sum_{k=1}^\infty c^k \, \bbE \, \big\langle e^d , 
	\big( \go_1 \cdots \go_k \big) \, e^d \big\rangle\,.
\end{equation} 
Notice that the non disordered case is recovered by setting
$\go_i \equiv I$, so that the diffusion constant
becomes $1/d + 2c/ (d(1-c))$. 
Assume now that $c > 0$ and let us switch the disorder on:
if we exclude the trivial case
when $\bbP \big( \go_1 e^d = e^d \big) = 1$,
we see that {\sl the diffusion constant decreases},
whatever the disorder law is.

We point out that by Schur's Lemma the relation $\overline r = cI$ is
fulfilled when the law of $r_1$ is {\sl conjugation invariant},
{\sl i.e.}, $\bP(r_1 \in \cdot) = \bP(h\, r_1\, h^{-1} \in \cdot)$
for every $h \in SO(d)$.

\item \rule{0pt}{1.15em}When the variables $\go_n$ are independent (and identically
distributed), and with no extra-assumption on $\overline r$,
by setting $\overline \go := \bbE(\go_1)$ we can write
\begin{equation}
	\gs^2 \;=\; \frac 1d \,+\, \frac 2d \, \sum_{k=1}^\infty
	\big\langle e^d , ( \overline{\go} \, \overline{r} )^k \, e^d \big\rangle
	\;=\; \frac 1d \,+\, \frac 2d \, 
	\bigg\langle e^d , \frac{\overline{\go} \, \overline{r}}
	{1 - \overline{\go} \, \overline{r}} \, e^d \bigg\rangle\,.
\end{equation}
Notice in fact that Assumption~\ref{ass:noise}
yields $\| \overline{r} \|_{\text{op}} < 1$, where
$\| \cdot \|_{\text{op}}$ denotes the operator norm
(see Section~\ref{sec:preparation}), hence the geometric series
converges.

\end{itemize}

\smallskip
\noindent
In the general case, the expression for the variance is not explicit, but of course
it can be evaluated numerically.

\smallskip

In order to get some intuition on the model, in particular  on 
the role of the disorder and  why it leads to \eqref{eq:sigma-def},
we suggest to have a look at
Appendix~\ref{sec:2d}, where we work out the computation of the asymptotic
variance of $X_n^{v_0}$ in the two-dimensional case, where elementary
tools are available because $SO(2)$ is Abelian.
As a matter of fact, these elementary tools would allow to 
prove for $d=2$ all the results we present in this paper.
However, the higher dimensional setting is
much more subtle and in particular
the proof of Theorem~\ref{th:main} for $d>2$ requires
more sophisticated techniques: in Section~\ref{sec:preparation},
using tensor analysis, we prove that Theorem~\ref{th:main} follows from
Assumption~\ref{ass:disorder} plus
a general condition of exponential convergence of some operator
norms, {\sl cf.} Hypothesis~\ref{hyp:main} below, and we then show that this condition is
a consequence of Assumption~\ref{ass:noise}.

\medskip

\begin{rem}\rm \label{rem:hom}
In the {\sl homogeneous} case, i.e., when disorder is absent, 
our method yields a proof  of the result in Theorem~\ref{th:main}
under a generalized irreducibility condition that
is weaker than Assumption~\ref{ass:noise} (see Appendix~\ref{app:hom}).
This generalized condition is fulfilled in particular whenever
the support of $Q$ generates a dense subset in $SO(d)$.
We point out that this last requirement is exactly the assumption under
which Theorem~\ref{th:main} (in the homogeneous case) was 
proven in~\cite{cf:Gor,cf:Roynette}.
\end{rem}

\medskip


\subsection{On strong decay of correlations}

\label{sec:corrs}
The persistence length ({\sl cf.} caption of Figure~\ref{fig:1}) 
does characterize the loss of the initial direction, but
from a probabilistic standpoint this is not completely satisfactory,
since other information could be carried on much further along the chain.
For this reason, we study the mixing properties of the variables 
$v_i^\omega$ (see \eqref{eq:X})
and this leads to a novel correlation length, that guaranties decorrelation of 
arbitrary local observables.
As we will see, we have only a bound on this new correlation length
and we can establish such a result 
only for a resticted (but sensible) class of models.

In order to state the result, let us introduce the
$\gs$-field $\cF_{m,n}^\go := \gs
(v_i^\go:\, m \le i \le n)$
for $m\in \N$, $n \in \N \cup \{\infty\}$ and for fixed $\go$. 
Then the mixing index $\ga^\go (n)$ of the sequence $\{ v_i^\go\}_i$
is defined for $n\in\N$ by
\begin{equation}
    \ga^\go (n) \, :=\, \sup \big\{ |\bP(A\cap B) - \bP(A)\bP(B)|:\;
    A \in \cF_{1,m}^\go, \, B\in\cF_{m+n,\infty}^\go,\, m \in \N \big\}.
\end{equation}
We work under either one of the following two hypotheses:

\smallskip
\newcounter{Hcount}
\begin{list}{H-\arabic{Hcount}.}
   {\usecounter{Hcount}
    }
\item  \label{H:1} The law $Q$ of $r_1$ is conjugation invariant,
{\sl i.e.}, $\bP(r_1 \in \cdot) = \bP(h\, r_1\, h^{-1} \in \cdot)$
for every $h \in SO(d)$,
and for some $n_0$ the law $Q^{*n_0}$
of $(r_1 \cdots r_{n_0})$ has an $L^2$ density with respect to the Haar measure
on $SO(d)$ (see \S\ref{sec:notations}).

\item \label{H:2}\rule{0pt}{1.2em}The law $Q$ of $r_1$ 
has an $L^2$ density with respect to the Haar measure on $SO(d)$.
\end{list}
\smallskip

\noindent
Assumption H-\ref{H:1} is sensibly weaker than  H-\ref{H:2}
(of course on the conjugation invariant measure),
however requiring an $L^2$ density is  quite a reasonable assumptions for thermal fluctuations.
Then we have

\medskip
\begin{proposition}
\label{th:mixing}
Under assumptions H-1 or H-2
there exist two constants
$C \in (0,\infty)$ and  $h \in (0,1)$ such that   
$\alpha^\go (n) \, \le \,  C\,  h^{n}$ for every $n$ and every $\go$
\end{proposition}
\medskip

The proof of  
Proposition~\ref{th:mixing} relies on Fourier analysis on $SO(d)$: it is given in
Section~\ref{sec:decay}, where one can find also an explicit characterization of the constant $ h$
(see \eqref{eq:h}).


\medskip

\section{The invariance principle}

\label{sec:preparation}

In this section we prove the invariance principle 
in Theorem~\ref{th:main}, including the formula \eqref{eq:sigma-def}
for the diffusion constant, under some abstract condition,
see Hypothesis~\ref{hyp:main} below,
which is then shown to follow from Assumption~\ref{ass:noise}.
Throughout the section we set
\begin{equation} \label{eq:phi}
	\gp_{m,n}^\go \;:=\; \go_m \, r_m \, 	\go_{m+1} \, r_{m+1} \,
	\ldots \, \go_n \, r_n \,, \qquad \quad m \le n \,,
\end{equation}
so that $v_j^\go = R\, \gp_{1,j}^\go \, e^d$ (see \eqref{eq:X}).
We recall that $v_0 = R\,e^d$ is an arbitrary element of $S^{d-1}$,
with $R \in SO(d)$, and that $Q$ denotes the law of $r_1$.


\smallskip

\subsection{Tensor products and operator norms}
\label{sec:tensors}

Unless otherwise specified, in this section the vector spaces
are assumed to be real ({\sl i.e.}, $\R$ is the underlying field)
and to have finite-dimension. 
The tensor product  of two vector spaces V and W can be introduced for example
by considering first 
the Cartesian product $V\times W$ and the (infinite-dimensional) vector space 
$\overline {V\times W}$ for which
the elements of $V\times W$ are a basis. Then the tensor product 
$V\otimes W$ is defined as the quotient space of $\overline {V\times W}$ under the 
equivalence relations 
\begin{gather*}
(v_1+v_2)\times w \;\sim\; v_1\times w \;+\; v_2\times w \,, \qquad \ \
v\times (w_1+w_2) \;\sim\; v\times w_1 \;+\; v\times w_2 \,, \\
c (v\times w) \;\sim\;  (cv)\times w \;\sim\; v \times (cw)\,,
\end{gather*}
for $c\in \R$, $v_{(i)}\in V$ and $w_{(i)}\in W$.
The equivalence class of $v\times w$ is denoted by $v\otimes w$ 
and we have the properties
$(v_1+v_2) \otimes w= v_1\otimes w+ v_2\otimes w$,
$v\otimes (w_1+w_2)= v\otimes w_1+ v\otimes w_2$
and $c (v\otimes w)= (cv)\otimes w= v\otimes (cw)$.
Given a basis $\{v_i\}_{i=1, \ldots , n}$ of $V$
and a basis $\{w_i\}_{i=1, \ldots , m}$ of $W$, $\{ v_i \otimes w_j\}_{i,j}$
is a basis of $V\otimes W$, which is therefore of dimension $nm$.
We stress that not every vector in $V \otimes W$ is of the form
$v \otimes w$ for some $v \in V$, $w \in W$.

A more concrete construction of $V \otimes W$ is possible
in special cases, e.g., when $V = W = \cL(\R^d)$,
the vector space of linear operators on $\R^d$ (that will be
occasionally identified with the corresponding representative
matrices in the canonical basis). 
In fact $\cL(\R^d) \otimes \cL(\R^d)$ is isomorphic to
$\cL(\cL(\R^d))$,
the space of all linear operators on $\cL(\R^d)$, and this
identification will be used throughout the paper.
Let us be more explicit: given $g,h \in \cL(\R^d)$,
we can view $g \otimes h$ as the linear operator sending $m \in \cL(\R^d)$ to
\begin{equation} \label{eq:goh}
	(g \otimes h) (m) \;:=\; g \, m \, h^* \,, 
	\qquad \text{that is} \qquad
	[(g \otimes h) (m)]_{ij} \;:=\; \sum_{k,l = 1}^d
	g_{ik} \, h_{jl} \, m_{kl}\,,
\end{equation}
where $(h^*)_{ij} =  h_{ji}$ is the adjoint of $h$. 
We are going to use this construction especially for $g, h \in SO(d)$, which of course
is not a vector space, but can be viewed as a subset of $\cL(\R^d)$.
A useful property of this representation of $g \otimes h$
as an operator is that
\begin{equation} \label{eq:comp}
	(g_1 \otimes h_1)(g_2 \otimes h_2) \;=\; (g_1 g_2) \otimes (h_1 h_2)\,,
\end{equation}
which is readily checked from \eqref{eq:goh}.
Another crucial fact is the following one: given 
$s_1, s_2 \in \cL(\R^d)^*$, the bilinear form
$(g,h) \mapsto s_1(g) s_2(h)$ can be written as a {\sl linear form}
$s_1 \otimes s_2$ on the tensor space $\cL(\R^d) \otimes \cL(\R^d)$, 
defined on product states $g \otimes h$ by
\begin{equation} \label{eq:tensorlin}
	(s_1 \otimes s_2)(g \otimes h) \;:=\; s_1(g) s_2(h)
\end{equation}
and extended to the whole space by linearity.
This linearization procedure is the very reason for introducing
tensor spaces, as we are going to see below.

Let us recall the definition and properties
of some operator norms. Given a vector space $V$ endowed with
a scalar product $\langle \cdot, \cdot \rangle$ and an operator
$A \in \cL(V)$, we define
\begin{equation} \label{eq:ophs}
	\| A \|_{\text{op}} \;:=\; \sup_{v,w \in V \setminus \{0\}} \, 
	\frac{|\langle w, A v \rangle|}{\|v\| \|w\|} \;=\;
	\sup_{v \in V \setminus \{0\}} \, \frac{\| Av \|}{\| v \|}\,, \qquad \
	\| A \|_{\text{hs}} \;:=\; \sqrt{\text{Tr}(A^* A)}\,,
\end{equation}
where $\text{Tr}(A)$ is the trace of $A$ and $A^*$ is the adjoint operator
of $A$, defined by the identity $\langle w, A v \rangle = \langle A^* w, v \rangle$
for all $v,w \in V$. If we fix an orthonormal basis $\{e^i\}_{i=1, \ldots, n}$
of $V$ and we denote by $A_{ij}$ the matrix of $A$ in this basis, we can write
$\| A \|_{\text{hs}}^2 \;:=\; \sum_{i,j} |A_{ij}|^2$.
It is easily checked that for all operators $A, B \in \cL(V)$ we have
\begin{equation} \label{eq:propophs}
	\|AB \|_{\text{op}} \;\le\; \|A\|_{\text{op}}\, \|B\|_{\text{op}}\,,
	\qquad 
	\| A \|_{\text{op}} \;\le\; \|A\|_{\text{hs}} \,, \qquad
	\| A B \|_{\text{hs}} \;\le\; 
	\|A\|_{\text{op}} \, \|B\|_{\text{hs}} \,.
\end{equation}
In what follows, the space $\cL(\R^d)$ is always equipped with the scalar product
$\langle v, w \rangle_{\text{hs}} := \text{Tr}(v^* w) = \sum_{i,j} v_{ij} w_{ij}$.
We can then give some useful bound on the operator norm of $g \otimes h$
acting on $\cL(\R^d)$: by \eqref{eq:goh} and \eqref{eq:propophs}
\begin{equation}
\label{eq:opineq}
\|g \otimes h\|_{\text{op}} \;=\;
\sup_{v \in \cL(\R^d) \setminus \{0\}}
\frac{\| g\, v\, h^* \|_{\text{hs}}}{\|v\|_{\text{hs}}}
\;\le\; \sup_{v \in \cL(\R^d) \setminus \{0\}}
\frac{\| g \|_{\text{op}} \|v\, h^* \|_{\text{hs}}}{\|v\|_{\text{hs}}}
\;\le\; \|g\|_{\text{op}} \, \|h\|_{\text{op}} \,,
\end{equation}
where we have used that $\|h^*\|_{\text{op}} = \|h\|_{\text{op}}$.

Let us denote by $\Gamma$ the orthogonal projection
on the subspace of symmetric operators in $\cL(R^d)$, defined 
for $v \in \cL(\R^d)$ by
\begin{equation} \label{eq:symm0}
	\Gamma(v) \;:=\; \frac 12 \,\big( v + v^*)\,, \qquad
	\text{i.e. } \qquad \Gamma(v)_{ij} \;=\; \frac 12 \,
	\big( v_{ij} + v_{ji} \big) \,.
\end{equation}
Of course $\Gamma \in \cL(\cL(\R^d))$, and for
any linear operator $m \in \cL(\cL(\R^d))$ we denote by $\overline m$ 
its symmetrized version:
\begin{equation} \label{eq:symm1}
	\overline m \;:=\; \Gamma \, m \, \Gamma \,.
\end{equation}
Note that $\overline{g \otimes g} 
\,=\, (g \otimes g)\, \Gamma \,=\, \Gamma \,(g \otimes g)$,
for every $g \in \cL(\R^d)$.

Finally, consider $s \in \cL(\R^d)^*$ of the form
$s(g) = \langle v, g\, w \rangle$, where $v, w$
are vectors in $\R^d$ with $\| v \| = \|w \| = 1$.
For every linear operator $m\in \cL(\cL(\R^d))$ we have
\begin{equation} \label{eq:symm2}
	(s \otimes s)(m) \;=\; (s \otimes s)(\Gamma \,m) \;=\; (s \otimes s)(m \, \Gamma)
	 \;=\; (s \otimes s)(\overline m) \,,
\end{equation}
as one easily checks using coordinates, since
$(s \otimes s)(m) = \sum_{ijkl} v_i \, v_j \, m_{ij,kl} \,  w_k \, w_l$.
It is also easily seen that
\begin{equation} \label{eq:two}
|(s \otimes s)(m)| \; \le \; \| m\|_{\text{op}} \,.
\end{equation}
These relations are easily generalized to higher order tensor products: in particular
\begin{equation} \label{eq:k}
	s^{\otimes 4}(m)
	\;=\; s^{\otimes 4} \big( m \, (\Gamma \otimes \Gamma) \big) 
	\qquad \text{and} \qquad
	| s^{\otimes 4}(m) | \;\le\; \| m \|_{\text{op}}\,,
\end{equation}
for every $m \in \cL(\R^d)^{\otimes 4}$.


\smallskip

\subsection{An abstract condition}

We are ready to state a condition on $Q$ that
will allow us to
prove the invariance principle in Theorem~\ref{th:main}.

Let us consider $\bE \phi^\go_{m,n}$, which is an element of $\cL(\R^d)$
(we recall that $\phi^\go_{m,n}$ is defined in \eqref{eq:phi}). We need to assume
that, when $k$ is large, $\bE \phi^\go_{n,n+k}$
is exponentially close to the zero operator on $\R^d$, uniformly in~$n$.
We are also interested in the asymptotic behavior of
$\bE \big[ \phi^{\go}_{n,n+k} \otimes \phi^{\go}_{n,n+k} \big]$,
which by \eqref{eq:goh} is a linear operator on $\cL(\R^d)$: we need that,
when $k$ is large and uniformly in~$n$,
the symmetrized version
$\bE \big[\, \overline{\phi^{\go}_{n,n+k} \otimes \phi^{\go}_{n,n+k}} \, \big]$
of this operator, cf. \eqref{eq:symm1} and \eqref{eq:symm0}, is exponentially close 
to the linear operator $\Pi$ defined as 
the orthogonal projection on the one dimensional
linear subspace of $\cL(\R^d)$ spanned by the identity matrix $(I_d)_{i,j} = \gd_{i,j}$,
$1 \le i,j \le d$, that is
\begin{equation} \label{eq:Pi}
	\Pi(v) \;:=\; \frac 1d \, \text{Tr}(v) \, I_d\,, \qquad
	v \in \cL(\R^d)\,.
\end{equation}
The reason why the operator $\Pi$ should have this form will be clear in \S\ref{sec:validity}.
Let us now state more precisely the hypothesis we make on~$Q$.

\medskip

\begin{hypothesis}
\label{hyp:main}
The law $Q$ of $r_1$ is such that, for $\bbP$--almost every $\omega$, we have
\begin{equation}
 C(\go) \;:=\;
 \sup_{n \ge 1} \; \left\{ \, \sum_{k=0}^\infty \big\|
 \bE \big[ \phi^{\go}_{n,n+k} \big] \big\|_{\rm{op}} \;+\;
 \sum_{k=0}^\infty \big\|
 \bE \big[\, \overline{\phi^{\go}_{n,n+k} \otimes \phi^{\go}_{n,n+k}} \, \big]
 \,-\, \Pi \, \big\|_{\rm{op}} \,\right\} \;<\; \infty\,.
 \end{equation}
\end{hypothesis}

\smallskip
\noindent
The next paragraphs are devoted to showing that
Theorem~\ref{th:main} holds if we assume Hypothesis~\ref{hyp:main}
together with Assumption~\ref{ass:disorder}. We then show
in \S\ref{sec:validity} that Hypothesis~\ref{hyp:main} indeed
follows from Assumption~\ref{ass:noise}.


\smallskip

\subsection{The diffusion constant}

We start identifying the diffusion coefficient $\gs^2$,
given by equation \eqref{eq:sigma-def}.
For any $\R^d$-valued random variable $Z$ we denote by $\mathrm{Cov}(Z)$ 
its covariance matrix: $\mathrm{Cov}(Z)_{i,j}= \cov(Z^i,Z ^j)$.

\begin{proposition}
\label{th:lem1}
If Hypothesis~\ref{hyp:main} and Assumption~\ref{ass:disorder} hold,
then for $\bbP$--almost every $\go$ and for every $v_0 \in S^{d-1}$ we have that
\begin{equation}
\label{eq:lem1}
\lim_{n \to \infty} \, \frac 1n \,
\mathrm{Cov}_{\bP} \left( 
X_n^{v_0, \go}
\right)
\; =\; \sigma^2 I_d \,,
\end{equation}
where
\begin{equation} \label{eq:sigma2-def}
\sigma^2 = 
\frac 1d \left( 1+2 \sum_{k=1}^\infty
\bbE \, \bE \, \langle e^d, \varphi^\omega_{1,k} e^d \rangle
\right),
\end{equation}
the series in the right-hand side being convergent.
\end{proposition}

\medskip

\noindent
{\it Proof.} 
By a standard polarization argument it is enough to prove that
for any $v  \in S^{d-1}$
\begin{equation}
\label{eq:lem1-var}
\lim_{n \to \infty} \, \frac 1n \,
\var_{\bP} \left( 
\langle v, X_n^{v_0, \go}\rangle
\right)
\, =\, \sigma^2 \,, \qquad \text{$\bbP(\dd \go)$--a.s.}\,,
\end{equation}
because
\begin{equation}
\mathrm{cov}_{\bP} \left( 
\langle e^i, X_n^{v_0, \go}\rangle,
\langle e^j, X_n^{v_0, \go}\rangle
\right) = \var_{\bP} \left(\left\langle \frac{e^i+ e^j}{\sqrt{2}}, X_n^{v_0, \go}\right\rangle \right) - \var_{\bP} \left(\left\langle  \frac{e^i- e^j}{\sqrt{2}}, X_n^{v_0, \go}\right\rangle \right).
\end{equation}
We recall that $X_n^{v_0, \go} = v_0 + \sum_{k=1}^n R \, \phi^\go_{1,k} \, e^d$,
where we set $v_0 = R e^d$ for some $R \in SO(d)$. For notational simplicity, we
redefine $X_n^{v_0, \go} := X_n^{v_0, \go} - v_0$ for the rest of the proof
(notice that this is irrelevant for the purpose of proving \eqref{eq:lem1-var}).
Introducing the notation
\begin{equation} \label{eq:def_sv}
	s_v(g) \;:=\; \langle v, R \, g \, e^d \rangle \,, \qquad \text{for } g \in \cL(\R^d)\,,
\end{equation}
we have the simple estimate
\begin{equation}
\label{eq:fctrlm}
\big\vert \bE \langle v,X_n^{v_0, \go}\rangle \big\vert
\; = \; \left \vert  s_v \left(\sum_{k=1}^n  \bE \, \phi_{1,k}^\go \right) \right\vert
\; \le \; \sum_{k \in \N} \big\| \bE \, \phi_{1,k}^\go \big\|_{\text{op}}
\; < \; \infty ,
\end{equation}
by Hypothesis~\ref{hyp:main}. This
shows that, in order to establish \eqref{eq:lem1-var}, it is sufficient to consider
\begin{equation}
\label{eq:diagnondiag}
\bE\big[ \langle v,X_n^{v_0, \go}\rangle ^2\big] \; =\;
\sum_{k=1}^n \bE \left[ \big (s_{v} ( \gp_{1,k}^\go )\big)^2
\right] \;+\; 2 \, \sum_{k=1}^n
\sum_{l=1}^{n-k} \bE \left[
s_{v} \left( \gp_{1,l}^\go \right)
s_{v} \left( \gp_{1,l+k}^\go \right)
\right]. 
\end{equation}
By \eqref{eq:tensorlin} and \eqref{eq:symm2} we can write $(s_{v} ( \gp_{1,k}^\go )\big)^2
= (s_v \otimes s_v)\big( \overline{\gp_{1,k}^\go \otimes \gp_{1,k}^\go} \big)$,
and by \eqref{eq:two} together with Hypothesis~\ref{hyp:main}
we can rewrite the fist sum as
\begin{equation}
	\sum_{k=1}^n \bE \left[ \left (s_{v} \left( \gp_{1,k}^\go \right)\right)^2 \right] 
	\;=\; s_v^{\otimes 2}\left(\sum_{k=1}^n  \bE \big[\,
	\overline{\gp_{1,k}^\go \otimes  \gp_{1,k}^\go}\, \big] \right) \;=\;
	n\, s_v^{\otimes 2}(\Pi)  \;+\; O(1) \,.
\end{equation}
In the same spirit the control the off-diagonal terms. We first observe that
by \eqref{eq:comp}
\begin{equation}
	\varphi_{1,l}^\omega \otimes \varphi^\omega_{1,l+k}
	\;=\; \varphi_{1,l}^\omega \otimes ( \varphi^\omega_{1,l} \, \varphi^\omega_{l+1,l+k} )
	\;=\; \big( \varphi_{1,l}^\omega \otimes \varphi^\omega_{1,l} \big)
	(I_d \otimes \varphi^\omega_{l+1,l+k}) \,,
\end{equation}
where $I_d \in \cL(\R^d)$ is the identity operator. Then
by \eqref{eq:tensorlin} and \eqref{eq:symm2} we can write
\begin{equation}
	s_{v}( \gp_{1,l}^\go) \, s_{v}( \gp_{1,l+k}^\go)
	\,=\, s_v^{\otimes 2} \big( \Gamma \, \gp_{1,k}^\go \otimes \gp_{1,l+k}^\go \big)
	\,=\, s_v^{\otimes 2} \big( \big( 
	\overline{\varphi_{1,l}^\omega \otimes \varphi^\omega_{1,l}} \big)
	(I_d \otimes \varphi^\omega_{l+1,l+k}) \big) \,.
\end{equation}
By \eqref{eq:two}, \eqref{eq:opineq} and \eqref{eq:propophs} we then obtain
$\bE \big[ s_{v}( \gp_{1,l}^\go) \, s_{v}( \gp_{1,l+k}^\go) \big]
\le \big\| \bE \varphi^\omega_{l+1,l+k} \big\|_{\text{op}}$,
hence by Hypothesis~\ref{hyp:main} it is the clear that
 \begin{equation}
 \label{eq:off-d1}
 \lim_{m \to \infty} \limsup_{n \to \infty}
 \frac 1n 
 \sum_{k=m}^n
\sum_{l=1}^{n-k} \bE \big[ s_{v}( \gp_{1,l}^\go) \, s_{v}( \gp_{1,l+k}^\go) \big] \, =\, 0 \,.
 \end{equation}
This allows us to focus on studying the limit as $n \to \infty$ and for fixed $k$ of 
\begin{equation}
\label{eq:withH}
	\frac 1n \, \sum_{l=1}^{n-k} \bE \big[ s_{v}( \gp_{1,l}^\go) \, s_{v}( \gp_{1,l+k}^\go) \big] 
	\;=\; \frac 1n \, \sum_{l=1}^{n-k} s_v^{\otimes 2}
	\big( \bE \big[ \, \overline{\varphi_{1,l}^\omega \otimes \varphi^\omega_{1,l}} \, \big]
	(I_d \otimes \bE \varphi^\omega_{l+1,l+k} ) \big)  \,.
\end{equation}
In this expression we can replace 
$\bE \big[\, \overline{\varphi_{1,l}^\omega \otimes \varphi^\omega_{1,l}}\, \big]$
by its limit $\Pi$ by making a negligible error (of order $1/n$),
by Hypothesis~\ref{hyp:main}. Furthermore, by the Ergodic Theorem
\begin{equation}
\label{eq:witherg}
\lim_{n \to \infty} \;
\frac {1} n  \,
\sum_{l=1}^{n-k}  s_v^{\otimes 2} \big( \Pi \,
(I_d \otimes \bE \varphi^\omega_{l+1,l+k} ) \big) \; =\;
s_v^{\otimes 2} \big( \Pi \, (I_d \otimes \bbE  \bE \varphi^\omega_{1,k})  \big) \,,
\quad \text{$\bbP(\dd \go)$--a.s.}\,.
\end{equation}
We have therefore proven that $\bbP(\dd \go)$--a.s. 
\begin{equation}
\label{eq:pre-pre-lem1}
\lim_{n \to \infty} \, \frac 1n \,
\var_{\bP} \langle v, X_n^{v_0, \go}\rangle \;=\;
s_v^{\otimes 2} (\Pi) \;+\;  2 \sum_{k=1}^\infty
s_v^{\otimes 2} \big( \Pi \, (I_d \otimes \bbE  \bE \varphi^\omega_{1,k}) \big) \,.
\end{equation} 
Let us simplify this expression: by
\eqref{eq:Pi} the representative matrix of $\Pi$ 
is $\Pi_{ij, kl} = \frac 1d \, \gd_{ij} \, \gd_{kl}$ and
by \eqref{eq:def_sv} we can write $s_v(g) = \sum_m (R^*v)_m \, g_{md}$, hence
\begin{align} \label{eq:asd}
s_v^{\otimes 2} (\Pi)  \;=\;
\sum_{ij} (R^*v)_{i} \, (R^*v)_{j} \, \Pi_{ij, dd} \;=\;
\frac 1d \, \| R^*v\|^2 \;=\; \frac 1d\,,
\end{align}
because $R \in SO(d)$ and $v$ is a unit vector.
The second term in the right hand side of \eqref{eq:pre-pre-lem1} is
analogous: setting for simplicity $m := \bbE  \bE \varphi^\omega_{l+1,l+k}
\in \cL(\R^d)$, the matrix of the operator $\Pi(I_d \otimes m)$ is given by
\begin{equation}
	[\Pi(I_d \otimes m)]_{ij, kl} \;=\; \sum_{a,b = 1}^d
	\Pi_{ij, ab} \, (I_d \otimes m)_{ab, kl} \;=\;
	\frac 1d \sum_{a,b = 1}^d \gd_{ij} \, \gd_{ab} \, \gd_{ak}
	\, m_{bl} \;=\; \frac 1d \, \gd_{ij} \, m_{kl} \,,
\end{equation}
hence
\begin{equation} \label{eq:asd2}
s_v ^{\otimes 2}[ \Pi (I_d \otimes m) ]
\;=\; \sum_{ij} (R^*v)_{i} \, (R^*v)_{j} \, [\Pi(I_d \otimes m)]_{ij, dd}
\;=\; \frac 1d \, m_{dd}\,.
\end{equation}
Since $m_{dd} = \bbE  \bE \langle e^d, \varphi^\omega_{1,k} e^d \rangle$,
we have shown that the right hand side of \eqref{eq:pre-pre-lem1} coincides with
the formula \eqref{eq:sigma2-def} for~$\gs^2$ and therefore
equation \eqref{eq:lem1-var} is proven.\qed

\smallskip


\smallskip

\subsection{The invariance principle}
\label{sec:inv_pr}

Next we turn to the proof of the full invariance principle.
The main tool is a projection of the increments
of our process $\{X^{v_0, \go}_n\}_n$ 
on martingale increments, to which the Martingale Invariance Principle can be applied.

We start setting $\hat s(g) := R \, g \, e^d$, so that
\begin{equation} \label{eq:increments}
	\hat s (\phi_{1,n}^\go) \;=\; v_n^\go \;=\; X_n^{v_0, \go} - X_{n-1}^{v_0, \go} \,,
\end{equation}
{\sl cf.} \eqref{eq:X} and \eqref{eq:phi}.
Recalling the definition \eqref{eq:def_sv} of $s_v(g)$, we have
$\hat s(g) = \sum_{i=1}^d s_{e^i}(g) \, e^i$.
For $n=1,2\dots$ we introduce the $\R^d$-valued process 
\begin{equation} \label{eq:Y}
	Y_{n} \;:=\; \hat s (\varphi_{1,n}^\omega) \,-\, \bE\big[ \hat s(\varphi_{1,n}^\omega) \big] \,.
\end{equation}
We now show that, for $\bbP$--a.e.~$\go$,
\begin{equation} \label{eq:decorr}
	\sup_{n \ge 1} \left\{ \sum_{k=0}^\infty 
	\big\| \bE \big[ Y_{n+k} \big| \cF_{1,n}^\go \big] 
	\big\|_{L^\infty(\bP; \R^d)}
	\right\} \; < \; \infty \,,
\end{equation}
where we recall that
$\cF^\go_{m,n} := \gs\big( \phi_{1,i}^\go:\, m \le i \le n \big)$.
Observe that $\bE \big[ Y_{n+k} \big| \cF_{1,n}^\go \big] = 
\hat s \big( \bE \big[ \varphi_{1,n+k}^\omega \big|  \cF_{1,n}^\go \big] 
- \bE \big[ \varphi_{1,n+k}^\omega \big] \big)$ and we can write
\begin{equation}
\begin{split}
	\bE\big[ \varphi_{1,n+k}^\omega \big|  \cF_{1,n}^\go \big] 
	\,-\, \bE \big[ \varphi_{1,n+k}^\omega \big] & \;=\;
	\varphi_{1,n}^\omega\, \bE\big[ \varphi_{n+1,n+k}^\omega \big|  \cF_{1,n}^\go \big] 
	\,-\, \bE \big[ \varphi_{1,n}^\omega \big] 
	\, \bE \big[ \varphi_{n+1,n+k}^\omega \big] \\
	& \;=\; \big( \varphi_{1,n}^\omega - \bE \big[ \varphi_{1,n}^\omega\big] \big) 
	\, \bE\big[ \varphi_{n+1,n+k}^\omega \big] \,.
\end{split}
\end{equation}
Since $\big\| \varphi_{1,n}^\omega - \bE \big[ \varphi_{1,n}^\omega\big] \big\|_{\text{op}} \le 2$,
we have
\begin{equation}
	\big\| \bE \big[ Y_{n+k} \big| \cF_{1,n}^\go \big] 
	\big\|_{L^\infty(\bP; \R^d)} 
	\;\le\; 2 \, \big\| \bE\big[ \varphi_{n+1,n+k}^\omega \big] \big\|_{\text{op}}\,,
\end{equation}
hence \eqref{eq:decorr} follows from Hypothesis~\ref{hyp:main}.

We are now ready to prove the invariance principle.
It is actually more convenient to
redefine $B^{v_0, \go}_{N}(t)$, which was introduced
in \eqref{eq:BN}, as $\frac{1}{\sqrt N} \, X^{v_0, \omega}_{ \lfloor Nt\rfloor }$, where 
$\lfloor a \rfloor \in \N \cup \{0\}$ denotes the integer part of $a$.
In this way, $B^{v_0, \go}_{N}(\cdot)$ is a process with trajectories in the Skorohod space
$D([0, \infty))$ of c\`adl\`ag functions, which is more suitable in order to apply
the Martingale Invariance Principle. However, since the limit process $\gs B$
has continuous paths,
it is elementary to pass from convergence in distribution on $D([0, \infty))$
to convergence on $C([0, \infty))$, thus recovering the original statement of Theorem~\ref{th:main}.

\medskip
\begin{theorem}
\label{th:inv-princ}
If Hypothesis~\ref{hyp:main} and Assumption~\ref{ass:disorder} hold,
then $\bbP(\dd \go)$--a.s. and for every choice of $v_0$
the $\R^d$-valued process $B^{v_0,\omega}_{N}$ 
converges in distribution on $C([0, \infty))$ to $\gs B$, 
where $B$ is a standard $d$-dimensional Brownian motion
and $\gs^2$ is given by \eqref{eq:sigma-def}.
\end{theorem}

\smallskip

\noindent {\it Proof.}
Let us set for $n \ge 1$
\begin{equation} \label{eq:UnZn}
U_n \,:=\, \sum_{k=0}^\infty \bE \big[ Y_{n+k}  \big| \cF_{1,n-1}^\go \big]
\quad \text{and} \quad
Z_n \,:=\, \sum_{k=0}^\infty \big( \bE \big[ Y_{n+k} \big| \cF_{1,n}^\go \big] 
\,-\, \bE\big[Y_{n+k} \big| \cF_{1,n-1}^\go\big] \big)\,,
\end{equation}
where we agree that $\cF_{1,0}^\go$ is the trivial $\gs$-field.
Note that $U_n$ and $Z_n$ are well-defined, because by equation \eqref{eq:decorr}
the series in \eqref{eq:UnZn} converge in $L^\infty(\bP;\R^d)$,
for $\bbP$--a.e.~$\go$.
The basic observation is that $\bE[Z_n|\cF_{1,n-1}^\go] = 0$,
hence $Z_n$ is a {\sl martingale difference sequence}, {\sl i.e.},
the process $\{T_n\}_{n \ge 0}$ defined by
\begin{equation} \label{eq:Tn}
	T_0 \;:=\; 0 \,, \qquad T_n \;:=\; \sum_{i=1}^n Z_i\,,
\end{equation}
is a $\{\cF^\go_{1,n}\}_n$--martingale (taking values in $\R^d$).
Moreover we have by construction
\begin{equation} \label{eq:telescopic}
Y_n  \;=\; \bE \big[ Y_n \big| \cF_{1,n}^\go \big] \;=\; Z_n \;+\;
(U_n - U_{n+1})\,,
\end{equation}
that is $Y_n$ is just $Z_n$ plus a telescopic remainder.
Therefore the process $\{T_n\}_n$ is very close to the original process $\{X^{v_0,\go}_n\}_n$,
because the variables $Y_n$ are nothing but the centered increments of the process
$\{X^{v_0, \go}_n\}_n$, see \eqref{eq:increments} and \eqref{eq:Y}.

For this reason, we start proving the invariance principle for the
rescaled process $T^N = \{T^N(t)\}_{t \in [0,\infty)}$
defined by $T^N(t) := \frac{1}{\sqrt N} T_{\lfloor Nt \rfloor}$.
By the Martingale Invariance Principle 
in the form given by \cite[Corollary~3.24, Ch.~VIII]{cf:JS},
the $\R^d$-valued process $T^N$ converges in law to $\tilde \gs B$,
where $\tilde \gs > 0$ and $B$ denotes a standard $\R^d$-valued Brownian motion,
provided the following conditions are satisfied:
\begin{itemize}
\item[(i)] the (random) matrix 
$(V_n)_{i,j} = \sum_{k=1}^n \bE \big[\, \langle e^i,Z_k\rangle \langle e^j, Z_k\rangle \,\big|
\,\cF_{1,k-1}^\go \,\big]$, with $1 \le i,j \le d$, is such that
\begin{equation}
\label{eq:cond1}
\frac 1n \, V_n \;\xrightarrow{\,n\to\infty\,}\; \tilde \sigma^2\, I_d \qquad 
\text{in $\bP$--probability\,;}
\end{equation}
\item[(ii)] the following integrability condition holds:
\begin{equation}
\label{eq:cond2}
\frac 1n \, \sum_{k=1}^n \bE\big[ \, |Z_k|^2 \,;\, |Z_k| > \varepsilon \sqrt{n} \, \big]
\;\xrightarrow{\,n\to\infty\,}\; 0 \,.
\end{equation}
\end{itemize}
The second condition is trivial because the variables $Z_n$ are bounded,
$\bbP(\dd\go)$--a.s.. The first condition requires more work. 
We first show that $\mathrm{var}_{\bP} \big( \frac 1n \, (V_n)_{i,j} \big) \to 0$ 
as $n \to \infty$, for all $i,j=1,\dots,d$ and for $\bbP$--a.e.~$\omega$,
and then we prove the convergence of $\bE \big[ \frac 1n V_n \big]$.

We start controlling the variance of $V_n$. By definition
$(V_n)_{i,j} \le \frac 12 ((V_n)_{i,i} + (V_n)_{j,j})$,
hence it suffices to show that
$\mathrm{var}_{\bP} \big( \frac 1n \, (V_n)_{i,i} \big) \to 0$
for every $i = 1, \ldots, d$.
We observet that $Z_n$ has a nice explicit formula:
\begin{equation}
Z_n = \hat s \left( \, \phi_{1,n-1}^\go
\, \big( \phi_{n,n}^\go - \bE [\phi_{n,n}^\go] \big) 
\left( \sum_{k=0}^\infty  \bE \big[ \phi_{n+1,n+k}^\go \big] \right) \right)\,,
\end{equation}
where we agree that $\phi_{n+1,n}^\go$ is the identity operator on $\R^d$
(this convention will be used throughout the proof).
Since $\hat s(g) = \sum_{i=1}^d s_{e^i}(g) \, e^i$,
where $s_v(g)$ is defined in \eqref{eq:def_sv}, a simple computation
then yields
\begin{equation} \label{eq:nice}
	\bE\big[\, \langle e^i,Z_n\rangle^2 \,\big|\, \cF_{1,n-1}^\go \,\big] 
	=\; s_{e^i}^{\otimes 2} \big(
	\bE\big[\, Z_n \otimes Z_n \,\big|\, \cF_{1,n-1}^\go \,\big] \big)
	\;=\; s_{e^i}^{\otimes 2}
	\left((\phi_{1,n-1}^\go)^{\otimes 2} \, \Theta^\go_n  \right) \,,
\end{equation}
where we have applied \eqref{eq:tensorlin} and \eqref{eq:comp} and we have set
\begin{equation}
\Theta^\go_n \,: =\,  \big( \bE[ \phi_{n,n}^\go \otimes \phi_{n,n}^\go] 
  \,-\, \bE [\phi_{n,n}^\go] \otimes \bE [\phi_{n,n}^\go] \big) \left(
  \sum_{k=0}^\infty \bE[\phi_{n+1,n+k}^\go] \otimes
  \sum_{l=0}^\infty  \bE [\phi_{n+1,n+l}^\go]  \right) \,.
\end{equation}
Applying \eqref{eq:nice} together with \eqref{eq:tensorlin} and \eqref{eq:comp}
we obtain
\begin{align}
	\mathrm{var}_{\bP} \bigg( \frac{(V_n)_{i,i}}{n} \bigg)
	& \;=\; \frac{1}{n^2} \, \bE \Bigg[ \Bigg( \sum_{1 \le k \le n} 
	s_{e^i}^{\otimes 2}
	\Big( \big( (\phi_{1,k-1}^\go)^{\otimes 2} 
	- \bE \big[ (\phi_{1,k-1}^\go)^{\otimes 2} \big] \big) \Theta^\go_k  \Big)
	\Bigg)^2\, \Bigg] \nonumber \\
	&\;\le\, \frac{2}{n^2} \, \sum_{1 \le k \le l \le n} 
	s_{e^i}^{\otimes 4} 	\Big( \bE \Big[ 
	\big( (\phi_{1,k-1}^\go)^{\otimes 2} - \bE \big[ (\phi_{1,k-1}^\go)^{\otimes 2} \big] 
	\big) \label{eq:bestia} \\
	& \qquad \qquad \qquad \qquad \qquad \quad
	\otimes \big( (\phi_{1,l-1}^\go)^{\otimes 2} - \bE\big[ (\phi_{1,l-1}^\go)^{\otimes 2}
	\big] \big) \Big]
	\big( \Theta^\go_l \otimes \Theta^\go_k \big) \Big) \,. \nonumber
\end{align}
Observe that by \eqref{eq:comp} we can write
\begin{equation}
\begin{split}
	(\phi_{1,l-1}^\go)^{\otimes 2} \;-\; \bE\big[ (\phi_{1,l-1}^\go)^{\otimes 2} \big]
	& \;=\; (\phi_{1,k-1}^\go)^{\otimes 2} \, (\phi_{k, l-1}^\go)^{\otimes 2}
	\;-\; \bE\big[ (\phi_{1,k-1}^\go)^{\otimes 2} \big] \,
	\bE\big[ (\phi_{k,l-1}^\go)^{\otimes 2} \big]\\
	& \;=\; \Big( (\phi_{1,k-1}^\go)^{\otimes 2} 
	- \bE\big[ (\phi_{1,k-1}^\go)^{\otimes 2} \big] \Big) \,
	\bE\big[ (\phi_{k,l-1}^\go)^{\otimes 2} \big] \\
	& \qquad \qquad \qquad + \; (\phi_{1,k-1}^\go)^{\otimes 2} \,
	\Big\{ (\phi_{k, l-1}^\go)^{\otimes 2} - \bE\big[ (\phi_{k,l-1}^\go)^{\otimes 2} \big] \Big\} \,,
\end{split}
\end{equation}
and notice that the term inside the curly brackets is independent
of $\cF^\go_{1,k-1}$ and vanishes when we take
the expectation. Therefore we have
\begin{equation}
\begin{split}
	& \bE \Big[ 
	\big( (\phi_{1,k-1}^\go)^{\otimes 2} - \bE \big[ (\phi_{1,k-1}^\go)^{\otimes 2} \big] 
	\big) \otimes \big( (\phi_{1,l-1}^\go)^{\otimes 2} - \bE\big[ (\phi_{1,l-1}^\go)^{\otimes 2}
	\big] \big) \Big] \\
	& \qquad \;=\;
	\bE \Big[(\phi_{1,k-1}^\go)^{\otimes 4} - 
	\bE \big[ (\phi_{1,k-1}^\go)^{\otimes 2} \big] ^{\otimes 2} \Big] \,
	\big( I \otimes \bE\big[ (\phi_{k,l-1}^\go)^{\otimes 2} \big] \big) \,,
\end{split}
\end{equation}
where we have applied again \eqref{eq:comp} and where $I$ denotes the
identity operator on $\cL(\R^d)$. We can therefore rewrite
the term in the sum in \eqref{eq:bestia} as
\begin{equation} \label{eq:combine}
\begin{split}
	& s_{e^i}^{\otimes 4} 	\Big( 
	\bE \Big[(\phi_{1,k-1}^\go)^{\otimes 4} - 
	\bE \big[ (\phi_{1,k-1}^\go)^{\otimes 2} \big]^{\otimes 2} \Big]
	\big( I \otimes \bE\big[ \phi_{k,l-1}^\go \otimes \phi_{k,l-1}^\go \big] \big)
	\big( \Theta^\go_l \otimes \Theta^\go_k \big) \Big) \\
	& =\, s_{e^i}^{\otimes 4} 	\Big( 
	\bE \Big[(\phi_{1,k-1}^\go)^{\otimes 4} - 
	\bE \big[ (\phi_{1,k-1}^\go)^{\otimes 2} \big]^{\otimes 2} \Big]
	\big( I \otimes \bE\big[\, \overline{\phi_{k,l-1}^\go \otimes \phi_{k,l-1}^\go} \, \big] \big)
	\big(\, \overline{\Theta^\go_l} \otimes \overline{\Theta^\go_k}\, \big) \Big) \,,
\end{split}
\end{equation}
where we have applied the first relation in \eqref{eq:k} together with
the following relations:
\begin{equation}
	\Theta^\go_l \, \Gamma \;=\; \overline{\Theta^\go_l} \qquad \text{and} \qquad
	\bE\big[ \phi_{k,l-1}^\go \otimes \phi_{k,l-1}^\go \big] \, \Theta^\go_k \, \Gamma
	\;=\; \bE\big[\, \overline{\phi_{k,l-1}^\go \otimes \phi_{k,l-1}^\go} \, \big] \, 
	\overline{\Theta^\go_k} \,,
\end{equation}
which follow from the fact that $(g \otimes g) \, \Gamma = \overline{g \otimes g}$
for every $g \in \cL(\R^d)$.

We know from Hypothesis~\ref{hyp:main} that when $l \gg k$ the operator
$\bE\big[\, \overline{\phi_{k,l-1}^\go \otimes \phi_{k,l-1}^\go} \, \big]$ is close to $\Pi$.
Furthermore, if we replace $\bE\big[\, \overline{\phi_{k,l-1}^\go \otimes \phi_{k,l-1}^\go} \, \big]$
by $\Pi$ inside \eqref{eq:combine} we get zero: in fact, since trivially
$g^{\otimes 2} \Pi = \Pi$ for every $g \in SO(d)$, we have
\begin{equation}
\begin{split}
	&\bE \Big[(\phi_{1,k-1}^\go)^{\otimes 4} \, - \, 
	\bE \big[ (\phi_{1,k-1}^\go)^{\otimes 2} \big]^{\otimes 2} \Big] 
	\big( I \otimes \Pi \big)\\
	& \;=\; \bE \big[(\phi_{1,k-1}^\go)^{\otimes 2} \otimes
	\big( (\phi_{1,k-1}^\go)^{\otimes 2}\, \Pi \big) \big] \; - \;
	\bE \big[ (\phi_{1,k-1}^\go)^{\otimes 2} \big] \otimes
	\bE \big[ (\phi_{1,k-1}^\go)^{\otimes 2} \, \Pi \big] \; =\; 0 \,.
\end{split}
\end{equation}
So it remains to take into account the contribution of the error 
$\bE\big[ \, \overline{\phi_{k,l-1}^\go \otimes \phi_{k,l-1}^\go} \, \big] - \Pi$
inside \eqref{eq:combine}.
However, using Hypothesis~\ref{hyp:main}, \eqref{eq:opineq} and the triangle inequality,
we have
\begin{equation}
	\left\|
	\bE \Big[(\phi_{1,k-1}^\go)^{\otimes 4} - 
	\bE \big[ (\phi_{1,k-1}^\go)^{\otimes 2} \big]^{\otimes 2} \Big]
	\right\|_{\text{op}} \, \le \, 2 \,, \qquad
	\left\|\Theta^\go_l \otimes \Theta^\go_k\right\|_{\text{op}}
	\, \le \, 4 \big(1+C(\go)\big)^4 \,,
\end{equation}
hence, using the second relation
in \eqref{eq:k}, from \eqref{eq:bestia} and \eqref{eq:combine} we obtain
\begin{equation}
\begin{split}
	\mathrm{var}_{\bP} \bigg( \frac{(V_n)_{i,i}}{n} \bigg)
	& \;\le\; \frac{8 \big(1+C(\go)\big)^4}{n^2} \, 
	\sum_{1 \le k \le n , \, m \ge 0} \big\|
	\bE\big[ \, \overline{\phi_{k,(k-1)+m}^\go \otimes \phi_{k,(k-1)+m}^\go} \, \big] - \Pi
	\big\|_{\text{op}}\\
	& \;\le\; \frac{8 \big(1+C(\go)\big)^5}{n} \,,
\end{split}
\end{equation}
having applied Hypothesis~\ref{hyp:main} again.
We have therefore shown that
$\mathrm{var}_{\bP} \big( \frac 1n \, (V_n)_{i,j} \big) \to 0$ as $n \to\infty$,
for all $1 \le i,j \le d$ and for $\bbP$--almost every~$\go$.

It remains to prove that $\bE \big[ \frac 1n \, V_n \big] \to \tilde\gs^2 \, I_d$
as $n\to\infty$ and to identify $\tilde\gs^2$.
Let us first note that by \eqref{eq:increments} and \eqref{eq:Y}
\begin{equation} \label{eq:babau}
	\sum_{k=1}^n Y_k \;=\; X_n^{v_0, \go} \;-\; \bE \big[ X_n^{v_0, \go} \big]
	\;=:\; \widetilde X_n
	\,.
\end{equation}
We also set $Z_n^i := \langle e^i, Z_n \rangle$,
$\widetilde X_n^i := \langle e^i, \widetilde X_n \rangle$
and $U_n^i := \langle e^i, U_n \rangle$ for short.
Since $\bE \big[ Z_n \big| \cF^\go_{1,n-1} \big] = 0$
and in view of \eqref{eq:telescopic}, we can write
\begin{equation} \label{eq:hence}
	\bE \big[ (V_n)_{i,j} \big] \;=\; \sum_{k=1}^n
	\bE \big[ Z_k^i \, Z_k^j \big] \;=\; \sum_{k,l=1}^n
	\bE \big[ Z_k^i \, Z_l^j \big] \;=\;
	\bE \big[ \big(\tilde X_n^i + U_{n+1}^i\big)\, \big(\tilde X_n^j + U_{n+1}^j \big) \big]
\end{equation}
(note that $U_1 = 0$). We recall that by Proposition~\ref{th:lem1} we have as $n \to \infty$,
for $\bbP$--a.e.~$\go$,
\begin{equation}
	\bE \big[\, \tilde X_n^i \, \tilde X_n^j \, \big] \;=\;
	\big( \text{Cov}_{\bP}( X_n^{v_0, \go}) \big)_{i,j} \;=\; n \, \gs^2 \, \gd_{i,j}
	\,+\, o(n) \,,
\end{equation}
where $\gs^2$ is given by \eqref{eq:sigma2-def} (equivalently by \eqref{eq:sigma-def}).
Since $\sup_n \|U_n\|_{L^\infty(\bP;\R^d)}<\infty$ by \eqref{eq:decorr},
it follows from \eqref{eq:hence} that as $n \to \infty$, for $\bbP$--a.e.~$\go$, we have
\begin{equation}
	\bE \big[ (V_n)_{i,j} \big] \;=\; \bE \big[\, \tilde X_n^i \, \tilde X_n^j \, \big]
	\,+\, o(n) \;=\; n \, \gs^2 \, \gd_{i,j} \,+\, o(n) \,.
\end{equation}
This completes the proof that the rescaled process $T^N = \{T^N(t)\}_{t \in [0,\infty)}$
converges in distribution as $N \to \infty$ to $\gs^2 \, B$, where
$\gs^2$ is given by \eqref{eq:sigma-def}.

It finally remains to obtain the same statement for
$B^{v_0, \go}_{N}(t) \,:=\, \frac{1}{\sqrt N} \, X^{v_0, \omega}_{ \lfloor Nt\rfloor}$.
Notice that by \eqref{eq:Tn}, \eqref{eq:telescopic} and \eqref{eq:babau} we can write
\begin{equation}
	\sup_{1 \le k \le n} \big\| X^{v_0, \omega}_k \,-\, T_k \big\| \;\le\;
	\sup_{1 \le k \le n} \big\| \bE \big[ X^{v_0, \omega}_k \big] \big\| \;+\;
	\sup_{1 \le k \le n} \big\| U_{k+1} \big\| \,,
\end{equation}
where $\|\cdot\|$ denotes the Euclidean norm in $\R^d$.
However the right hand side is bounded in $n$ in $L^\infty(\bP;\R^d)$,
for $\bbP$--a.e.~$\go$  (for the first term
see \eqref{eq:fctrlm} while for the second term
we already know that $\sup_n \|U_n\|_{L^\infty(\bP;\R^d)}<\infty$).
Therefore $\sup_{t \in [0,M]} \| B^{v_0, \go}_{N}(t) - T^N(t) \|
\le (const.)/\sqrt{N}$ for every $M > 0$, and the proof is completed.\qed


\smallskip

\subsection{Proof of Theorem~\ref{th:main}}

\label{sec:validity}

We now show that the abstract condition expressed by Hypothesis~\ref{hyp:main}
is a consequence of Assumption~\ref{ass:noise}. In view of Theorem~\ref{th:inv-princ},
this completes the proof of Theorem~\ref{th:main}.

\smallskip

We start by controlling $\bE \big[ \phi^\go_{n,n+k} \big]$, which is quite easy:
the independence of the $r_i$ yields
\begin{equation} \label{eq:val1}
	\bE \big[ \phi^\go_{n,n+k} \big] \;=\;
	\go_n \, \bE(r_1) \, \go_{n+1} \, \bE(r_1) \, \cdots
	\go_{n+k} \, \bE(r_1)\,\,.
\end{equation}
It is clear that $\| \bE(r_1) \|_{\text{op}} \le 1$.
We now show that Assumption~\ref{ass:noise} yields $\| \bE(r_1) \|_{\text{op}} < 1$,
so that for every $\go$ we have
\begin{equation} \label{eq:analogy}
	\sum_{k=1}^\infty \big\| \bE [\phi^\go_{n,n+k}] \big\|_{\text{op}}
	\;\le\; \sum_{k=1}^\infty \big\| \bE(r_1)^k \big\|_{\text{op}}
	\;\le\; \sum_{k=1}^\infty \| \bE(r_1)\|_{\text{op}}^k \;<\; \infty \,.
\end{equation}
To prove that $\| \bE(r_1) \|_{\text{op}} < 1$, we argue by contradiction:
if $\| \bE(r_1) \|_{\text{op}} = 1$ there would exist two vectors
$x,y \in S^{d-1}$ such that
\begin{equation}
	1 \;=\; \langle y, \bE(r_1) \, x \rangle \;=\; \int_{SO(d)} \langle y, g \, x \rangle
	\, Q(\dd g) \,.
\end{equation}
Since $|\langle y, g \, x \rangle | \le 1$, for this equality to hold it is
necessary that $g \, x = y$ for $Q$--almost every~$g \in SO(d)$. 
Setting $V := \{\gl\, x:\, \gl \in \R\}$ and $W := \{\gl\, y: \, \gl \in \R\}$,
this would mean that $g \, V = W$ for $Q$--almost every~$g \in SO(d)$, which is
in contradiction with Assumption~\ref{ass:noise}.

\smallskip

Next we turn to the analysis of
$\bE \big[\, \overline{\phi^\go_{n,n+k} \otimes \phi^\go_{n,n+k}} \,\big]$,
which is a linear operator on the vector space $\cL(\R^d)$,
equipped with the standard scalar product
$\langle v, w \rangle_{\text{hs}} = \text{Tr}(v^* w)$.
We decompose $\cL(\R^d) = H_1 \oplus H_s^0 \oplus H_a$ as a sum of
the orthogonal subspaces consisting respectively of the multiples of the identity,
of the symmetric matrices with zero trace and of the antisymmetric matrices:
\begin{gather*}
	H_1 \;:=\; \big\{ \gl\, I_d : \ \gl \in \R \big\} \,, \qquad \quad
	H_s^0 \;:=\; \big\{ v \in \cL(\R^d): \ v^* = v \ \text{and} \ \text{Tr}(v) = 0 \big\}\,,\\
	H_a \;:=\; \big\{ v \in \cL(\R^d): \ v^* = -v \big\} \,.
\end{gather*}
All of these subspaces are invariant under $g \otimes g$, for every $g \in \cL(\R^d)$,
hence they are invariant under $\bE \big[\, \overline{\phi^\go_{n,n+k} \otimes \phi^\go_{n,n+k}} \,\big]$.
We recall that $\Pi$ is the orthogonal projection on $H_1$, cf. \eqref{eq:Pi},
while $\Gamma$ is the orthogonal projection on $H_1 \oplus H_s^0$, cf. \eqref{eq:symm0}.
Since $\Pi$ and $\bE \big[\, \overline{\phi^\go_{n,n+k} \otimes \phi^\go_{n,n+k}} \,\big]$
are zero on $H_a$ and they coincide on $H_1$,
$\big\| \bE \big[\, \overline{\phi^\go_{n,n+k} \otimes \phi^\go_{n,n+k}} \,\big]
\, -\, \Pi \big\|_{\text{op}}$ is nothing but the operator norm of
$\bE \big[\phi^\go_{n,n+k} \otimes \phi^\go_{n,n+k}\big]$
restricted to the subspace $H_s^0$, therefore with obvious notation we can write
for every~$\go$
\begin{equation}
	\sum_{k=1}^\infty \big\| \bE \big[\, 
	\overline{\phi^\go_{n,n+k} \otimes \phi^\go_{n,n+k}} \,\big] - 
	\Pi \big\|_{\text{op}} \;=\; \sum_{k=1}^\infty 
	\big\| \bE [\phi^\go_{n,n+k} \otimes \phi^\go_{n,n+k}] \big\|_{H_s^0, \text{op}} \,.
\end{equation}
However from \eqref{eq:comp} and from the fact that the $r_i$ are independent and identically
distributed we have
\begin{equation} \label{eq:multip}
	\bE \big[ \phi^\go_{n,n+k} \otimes \phi^\go_{n,n+k}\big] 
	\;=\; (\go_n \otimes \go_n) \, \bE[ r_1 \otimes r_1] \, \cdots \,
	(\go_{n+k} \otimes \go_{n+k}) \, \bE[ r_1 \otimes r_1] \,,
\end{equation}
hence
\begin{equation} \label{eq:analogy3}
	\sum_{k=1}^\infty \big\| \bE \big[\, 
	\overline{\phi^\go_{n,n+k} \otimes \phi^\go_{n,n+k}} \,\big] 
	- \Pi \big\|_{\text{op}}
	\;\le\; \sum_{k=1}^\infty 
	\big( \big\| \bE[ r_1 \otimes r_1] \big\|_{H_s^0, \text{op}}\big)^k \,.
\end{equation}

We are finally left with showing that $\| \bE [r_1 \otimes r_1]\|_{H_s^0, \text{op}} < 1$.
Let us assume by contradiction that there exist
$v,w \in H_s^0$ with $\|v\|_{hs} = \|w\|_{hs} = 1$ such that
\begin{equation}
	1 \;=\; \langle w, \bE [r_1 \otimes r_1] \, v \rangle_{\text{hs}} \;=\; 
	\int_{SO(d)} \langle w, g\, v\, g^* \rangle_{\text{hs}} \, Q(\dd g) \,.
\end{equation}
However $\| g\, v\, g^* \|_{\text{hs}} = \| v \|_{\text{hs}} = 1$,
hence $\langle w, g\, v\, g^* \rangle_{\text{hs}} \le 1$ and
we must have $w = g\, v\, g^* = g\, v\, g^{-1}$ for $Q$--a.e.~$g$ in $SO(d)$.
In particular, the matrices $v$ and $w$ are similar and therefore
they have the same eigenvalues $\gl_1, \ldots, \gl_k$, with $k \le d$.
Recall that by the spectral theorem $v$ and $w$ are diagonalizable.
Denoting by $K_v$ and $K_w$ respectively the eigenspaces of $v$ and $w$
corresponding to $\gl_1$, we have that $1 \le \dim (K_v) = \dim (K_w) \le d-1$,
where the last inequality follows from the fact that $v$ and $w$,
having zero-trace and not being identically zero, 
cannot be multiples of the identity.
Let us now fix $g$ such that $w = g\, v\, g^{-1}$ and
take an arbitrary $x \in g \, K_v$: since $g^{-1} x \in K_v$ we have
\begin{equation}
	w\,x \;=\; g \, v \, g^{-1} \, x \;=\; g \, ( \gl_1 \, g^{-1} \, x) \;=\; \gl_1 \, x \,,
\end{equation}
which yields $x \in K_w$. Therefore $g \, K_v \subseteq K_w$ and since
the two subspaces have the same dimension we must have $g \, K_v = K_w$,
for $Q$--almost every~$g$. This being in contradiction with Assumption~\ref{ass:noise},
we have indeed that $\| \bE [r_1 \otimes r_1]\|_{H_s^0, \text{op}} < 1$ and the
proof of Theorem~\ref{th:main} is completed.\qed


\medskip

\section{Decay of correlation}
\label{sec:decay}


\smallskip

\subsection{General notations}

\label{sec:notations}

We denote by $\gl$ the normalized Haar measure on $SO(d)$. We recall that $\gl$ is 
the only probability measure that is left- and right-invariant,
{\sl i.e.}, such that $\gl(A g) = \gl (gA) = \gl(A)$
for all $g \in SO(d)$ and (measurable) $A \subseteq SO(d)$.
In the special case $d=3$, $\gl$ describes a (random) rotation around the vector 
$w$ of angle $\theta$, where $w$ is uniform on $S^{2}$ and $\theta$ is uniform on $[0, 2\pi)$.
For more on the Haar measure we refer to \cite{cf:Halmos}.  

We recall that $Q$ denotes the law of $r_1$.
For fixed $\go$,
we denote by $L_{m,n}^\go$ the law of $\gp_{m,n}^\go$ under $\bP$, so that
for any bounded and measurable function $F: SO(d) \to \R$
\begin{equation}
\bE\big[ F(\gp_{m,n}^\go) \big] \;=\; \int_{SO(d)} F(g) \, L_{m,n}^\go(\dd g) \,.
\end{equation}
We also set
\begin{equation}\label{eq:Ek}
E^\go(k) \;:=\; 2 \, \sup_{n} \, \| L^\go_{n+1, n+k} - \gl \|_{TV} \,,
\end{equation}
where the total variation (TV) distance between
the probability measures $\mu$ and $\nu$  
is defined as $\Vert \mu -\nu \Vert_{\text{TV}} :=
\sup_A \vert \mu(A) -\nu (A)\vert$. We observe that
$\Vert \mu -\nu \Vert_{\text{TV}}$ coincides with
$\frac 12 \sup_{\vert g\vert \le 1}\int g \dd \mu - \int g \dd \nu$,
in particular if $\mu$ is absolutely continuous with respect to $\nu$,
with $f := \dd \mu/\dd \nu$, we have $\| \mu - \nu \|_{\text{TV}} =
\frac 12 \int |f - 1 | \, \dd \nu$.


\subsection{Reminders of harmonic analysis on compact groups}

Throughout this section, we assume that $G$ is a compact topological group,
equipped with the Borel $\gs$-field,
and $\lambda$ is the normalized Haar measure on $G$
(of course we have in mind the specific case where $G=SO(d)$, $d \ge 2$).
We start recalling some basic facts about harmonic analysis on $G$,
taking inspiration from \cite{cf:Hewitt1,cf:Hewitt2}.

Given a (complex) Hilbert space $H$, a representation of $G$ on $H$ is a 
group homomorphism $U:G \to \cB(H)$, {\sl i.e.}, $U(gh) = U(g)U(h)$ for all $g,h \in G$,
where $\cB(H)$ denotes the set of bounded linear operators
from $H$ to itself. The representation $U$ is said to be:
\begin{itemize}
\item {\sl continuous} if the map $g \mapsto \langle x, U(g) y \rangle$
from $G$ to $\C$ is continuous, for all $x,y \in H$;
\item {\sl irreducible} if there is no closed subspace $M$ of $H$ such that 
$U(g) M \subseteq M$ for every $g \in G$, 
except the trivial case when $M=\{0\}$ or $M=H$;
\item {\sl unitary} if $U(g)$ is a unitary operator for every $g \in G$, {\sl i.e.},
$\langle U(g) x, U(g) y \rangle = \langle x,y \rangle$ for all $x,y \in H$, where
$\langle\cdot,\cdot\rangle$ denotes the scalar product in $H$ (that we
take skew-linear in the first argument and linear in the second).
\end{itemize}
Finally, two representations $U$, $U'$ of $G$ on the Hilbert spaces $H$, $H'$
are said to be {\sl equivalent} if there exists a linear isometry $T: H \to H'$
such that $U(g) = T^{-1} U'(g) T$ for every $g \in G$.
The set of equivalence classes of continuous, irreducible, unitary representations
of $G$ is denoted by $\Sigma$, which is a countable
set (sometimes called the {\sl dual object} of $G$).

We point out that, since $G$ is compact, all irreducible representations are
finite dimensional, that is, they act on a finite dimensional Hilbert space.
Given $\alpha \in \Sigma$, we denote by $U^\alpha$
an arbitrary representation in the class $\alpha$, acting on the Hilbert space $H_\alpha$
of finite dimension $d_\alpha \in \N$. 
In each space $H_\alpha$ we fix an (arbitrary) orthonormal basis
$\{\zeta^\alpha_i, i=1,\dots,d_\alpha\}$ and we denote by
$u^{\alpha}_{ij}(g) = \langle \zeta^\alpha_i, U^\alpha(g) \zeta_j^\alpha \rangle$
the matrix of $U^\alpha(g)$ on this basis.
Notice that $u^{\alpha}_{ij}(\cdot)$ is a continuous function from $G$ to $\C$.
We have the following orthogonality relations, valid for all
$\ga,\beta \in \Sigma$, $1 \le j,k \le d_{\alpha}$, $1 \le l,m \le d_{\beta}$:
\begin{equation}
  \label{eq:ortho-coeffs}
\int_G \overline{u^\alpha_{jk}(g)} \, u^\beta_{lm}(g) \, \lambda(\dd g) \;=\;
\frac {1}{d_\alpha} \, \delta_{\ga \gb} \, \delta_{jl} \, \delta_{km} \,,
\end{equation}
where $\overline{x}$ denotes the complex conjugate of $x$ and $\delta_{ij}$
is the Kronecker delta.
Therefore $\{\sqrt{d_\alpha} \, u^\alpha_{i,j}(\cdot)\}_{\alpha \in \Sigma,\,
1 \le i,j \le d_\alpha}$ is an orthonormal set in $L^2(G, \dd\lambda)$.
A crucial result is that it is also complete, {\sl i.e.},
the functions $u^\ga_{i,j}(\cdot)$ span $L^2(G, \dd\lambda)$,
by the the Peter-Weil Theorem.

Next we introduce the {\sl Fourier transform} $\hat\mu$ of a
probability measure $\mu$ on $G$, which is the element of the space
${\sf S} := \prod_{\alpha \in \Sigma} \cB(H^\alpha)$ defined by
\begin{equation} \label{eq:fourier}
	\hat \mu(\alpha) \;:=\; \int_G U^{\alpha}(g) \,\mu(\dd g) \,, \qquad
\alpha \in \Sigma\,.
\end{equation}
More explicitly, $\hat \mu(\alpha)$ is the linear operator
acting on $H^\alpha$ whose matrix
in the basis $\{\zeta^\alpha_i\}_{i}$ is given by
$\hat \mu(\alpha)_{i,j} = \int_G u^\alpha_{i,j}(g)\, \mu(\dd g)$,
for $\alpha \in \Sigma$ and $1 \le i,j \le d_\alpha$.

It follows directly from the definition \eqref{eq:fourier} and \eqref{eq:ophs}
that $\|\hat \mu(\alpha)\|_{\text{op}} \le 1$ for every
probability measure $\mu$ on $G$ and for every
$\alpha \in \Sigma$. As a matter of fact, when $G$ is connected, this
inequality is strict for a large class of $\mu$,
as we show in the following lemma (where we denote by
$\alpha = 0$ the trivial representation, with
$H_0 = \C$ and $U^0(g) = 1$ for every $g \in G$).

\smallskip
\begin{lemma} \label{lem:1}
Let $\mu$ be a probability measure on $G$ with support $V$.
Assume that $V^{-1} V := \{h^{-1}g:\, h,g \in V\}$ generates a dense set in $G$, 
{\sl i.e.}, the set $\union_{n=1}^\infty (V^{-1} V)^n$ is dense in $G$.
Then $\|\hat \mu(\alpha)\|_{\text{op}} < 1$ for
every $\alpha \in \Sigma$, $\alpha \ne 0$.
\end{lemma}
\smallskip

\proof
Suppose that $\|\hat\mu(\alpha)\|_{\text{op}} = 1$. Then there must
exist $x,y \in H^\alpha$ with $\|x\| = \|y\| = 1$ 
such that $ \langle y, \hat \mu(\alpha)\, x \rangle  = 1$.
Now
\begin{equation}
1= \Re \langle y , \hat \mu(\alpha) \,x \rangle 
= \int \Re \langle y,U^\alpha(g)\, x \rangle \nu(dg)
\end{equation}
The function $r(g) = \Re \langle y,U^\alpha(g)\,x\rangle$ is real and such that 
$r(g)\le 1$, and so must be constant on the support of $\mu$ and equal to $1$. 
This implies that $U^\alpha(g)\, x = y$ for any $g\in V$, hence
\begin{equation}
	U^\ga(h^{-1}g)\, x \;=\; U^\ga(h^{-1})\, U^\ga(g)\, x 
	\;=\; U^\ga(h^{-1})\, y \;=\; x
\end{equation}
for all $g,h \in V$. This means that the relation $U^\ga(g) x = x$
holds for all $g \in V^{-1}V$ and hence for all $g \in \union_{n=1}^\infty
(V^{-1} V)^n$. By assumption the latter set is dense in $G$ and the continuity
of the representation $U^\ga$ yields that $U^\ga(g) x = x$ for all $g \in G$,
which is impossible unless $\ga$ is the trivial representation.
\qed
\medskip

We conclude this paragraph noting that
the Fourier transform provides
an easy tool to check whether a probability measure
$\mu$ has an $L^2$ density with respect
to the Haar measure $\lambda$. More precisely, we have the following

\smallskip
\begin{lemma}[Fourier inversion theorem] \label{lem:Fourinv}
A probability measure $\mu$ on $G$ is such that
\begin{equation} \label{eq:condlemma}
	\sum_{\alpha \in \Sigma} d_\alpha \, \|\somb \mu(\alpha)\|_{\mathrm{hs}}^2 
	\;<\; \infty
\end{equation}
if and only if it is absolutely continuous with respect to $\lambda$ with density
in $L^2(G,\dd\lambda)$.
In this case, the density $f = \dd \mu/\dd \gl$ is given by
\begin{equation} \label{eq:repres}
	f(g)  \;=\; 
\sum_{\alpha\in\Sigma} d_\alpha \tr( \somb \mu(\alpha) \, U^\alpha(g)^* )
\;=\; \sum_{\alpha\in\Sigma} \sum_{1 \le i,j \le d_\ga}
d_\alpha \, \somb\mu(\alpha)_{i,j} \, 
\overline{u^\alpha_{i,j}(g)} \,,
\end{equation}
where the series converges in $L^2(G,\dd\gl)$.
\end{lemma}
\smallskip

\proof
Since $\{\sqrt{d_\alpha} \, u^\alpha_{i,j}(\cdot)\}_{\alpha \in \Sigma,\,
1 \le i,j \le d_\alpha}$ is a complete orthonormal set in $L^2(G, \dd\lambda)$,
the same is true if we replace
$u^\alpha_{i,j}(\cdot)$ by $\overline{u^\alpha_{i,j}(\cdot)}$,
therefore condition \eqref{eq:condlemma} guarantees that the right hand side
of \eqref{eq:repres} does define a function $f \in L^2(G,\dd \gl)$.
Consider then the (a priori complex) measure $\dd\nu := f \, \dd \gl$.
Using \eqref{eq:ortho-coeffs} it is easy to check that
\begin{equation}
	\somb \nu (\ga)_{i,j} \;:=\; \int_G u^\ga_{i,j}(g) \, \nu(\dd g)
	\;=\; \somb \mu(\ga)_{i,j}\,,
\end{equation}
for all $\ga \in \Sigma$ and $1 \le i,j \le d_\ga$.
By Theorem (27.42) of \cite{cf:Hewitt2} this implies that $\mu=\nu$.
Vice versa, if a function $f$ is in $L^2(G, \dd \gl)$,
the right hand side of \eqref{eq:repres} is nothing
but its Fourier series in the orthonormal set
$\{\sqrt{d_\alpha} \, \overline{u^\alpha_{i,j}(\cdot)}\}_{\alpha \in \Sigma,\,
1 \le i,j \le d_\alpha}$, hence relation \eqref{eq:condlemma} holds true.
Finally, the second equality in \eqref{eq:repres} is easily checked.
\qed
\smallskip


\smallskip

\subsection{Exponential decay of the total variation norm}

In this subsection we need to assume that $G$ is also connected
(which is of course the case for $G=SO(d)$).
We show that, assuming hypothesis H-\ref{H:1} or  hypothesis H-\ref{H:2}
({\sl cf.} \S~\ref{sec:corrs}), 
for $\bbP$--a.e. $\go$, we have
\begin{equation}
	\sum_{k \in \N} E^\go(k) \;<\; \infty\,,
\end{equation}
where we recall that $E^\go(k)$ has been introduced in \eqref{eq:Ek}.
As a matter of fact, we are going to prove the much stronger result
that there exist positive constants $c_1, c_2$ such that
\begin{equation} \label{eq:toprove}
	\sup_{\go} \, E^\go(k) \;\le\; c_1 \, e^{-c_2\, k}\,,
	\qquad \text{ for all } k \in \N \,.
\end{equation}

\smallskip

It is convenient to introduce the convolution $\mu*\nu$ of
two probability laws $\mu,\nu$ on $G$ by
\begin{equation}
	(\mu*\nu)(A) \;:=\; \int_G \mu(A h^{-1}) \, \nu(\dd h)
	\;=\; \int_G \nu(g^{-1}A) \, \mu(\dd g) \,,
\end{equation}
so that if $X$, $Y$ are two independent random elements of
$G$ with marginal laws $\mu$, $\nu$, then $\mu *\nu$ is the law of $XY$.
Therefore we can express $L^\go_{m,n}$ as
\begin{equation} \label{eq:multconv}
	L^\go_{m,n} \;=\; \gd_{\go_m} * Q * \gd_{\go_{m+1}} * Q * \ldots
	* \gd_{\go_n} * Q\,,
\end{equation}
where $\gd_{g}$ denotes the Dirac mass at $g \in G$.
We stress that in general the convolution is not commutative.
A basic property is that 
$\widehat{\mu*\nu}(\ga) = \somb\mu(\ga) \, \somb\nu(\ga)$ for every $\ga \in \Sigma$,
or more explicitly $\widehat{\mu*\nu}(\ga)_{i,j} = 
\sum_{k=1}^{d_\ga} \somb\mu(\ga)_{i,k} \, \somb\nu(\ga)_{k,j}$,
as one easily checks from \eqref{eq:fourier}.


In the next crucial lemma we give an explicit bound on $E^\go(k)$
in terms of the Fourier transform $\somb Q$ of $Q$.
We recall that we denote by $\ga = 0$ the trivial representation.

\smallskip
\begin{lemma} \label{lem:crucial}
The following relation holds true for every $k \in \N$:
\begin{equation} \label{eq:condE}
	\left( \sup_\go \, E^\go(k) \right)^2 
	\;\le\; \sum_{\alpha \in \Sigma, \, \alpha \ne 0}
	d_{\alpha} \,  \|\somb Q(\alpha)\|_{\mathrm{hs}}^{2} \,
	\|\somb Q(\alpha)\|_{\mathrm{op}}^{2(k-1)} \,.
\end{equation}
\end{lemma}
\smallskip

\proof
From \eqref{eq:multconv} we can write
\begin{equation} \label{eq:mydec}
\somb L^\go_{n+1,n+k}(\alpha) \;=\;
U^\alpha(\go_{n+1}) \, \somb Q(\alpha) \,\cdots \,
U^\alpha(\go_{n+k}) \, \somb Q(\alpha)\,,
\end{equation}
and using the inequalities in \eqref{eq:propophs} we get
\begin{equation}
\|\somb L^\go_{n+1,n+k}(\alpha)\|^2_{\mathrm{hs}} \;\le\;
\|\somb Q(\alpha) \|^{2}_{\mathrm{hs}} \,
\|\somb Q(\alpha) \|^{2(k-1)}_{\mathrm{op}} \,
\prod_{i=1}^k \| U^\alpha(\omega_{n+i} ) \|^2_{\mathrm{op}}
\;\le\; \|\somb Q(\alpha) \|^{2}_{\mathrm{hs}} \,
\|\somb Q(\alpha) \|^{2(k-1)}_{\mathrm{op}} \,,
\end{equation}
where we used that $\| U^\alpha(\go_{n+i} )\|^2_{\mathrm{op}} = 1$
because the representation is unitary.
Now assume that the right hand side of \eqref{eq:condE} is finite
(otherwise there is nothing to prove). By Lemma~\ref{lem:Fourinv}, 
$L^\go_{n+1,n+k}$ has a density 
$f \in L^2(G, \dd\gl)$ with respect to $\lambda$, therefore
by Jensen's inequality we can write
\begin{equation}
	4 \| L^\go_{n+1,n+k} - \lambda \|^2_{TV} \;=\; 
	\left(\int_G |f-1| d\lambda\right)^2 \;\le\;
	\int_G f^2 \, \dd \gl \; - \; 1  \;=\;
	\sum_{\alpha \in \Sigma, \, \ga \ne 0} d_{\alpha} \,
	\|\somb L^\go_{n+1,n+k} (\alpha) \|^2_{\mathrm{hs}} \,,  
\end{equation}
where in the last equality we have used Parseval's identity,
observing that $\langle f, u^\ga_{i,j} \rangle =
\big(\somb L^\go_{n+1,n+k} (\alpha)\big)_{i,j}$ and
that trivially $\somb \mu(0) = 1$ for every probability measure $\mu$
on $G$. Recalling the definition \eqref{eq:Ek} of $E^\go(k)$,
relation \eqref{eq:condE} is proven.\qed
\medskip

\noindent
{\it Proof of \eqref{eq:toprove} under
Hypothesis H-\ref{H:2}.} 
Let us set $f := \dd Q/\dd \gl \in L^2(G, \dd\gl)$.
By Parseval's identity we have
\begin{equation}
\| f \|^2_2 \;:=\; \int_G f^2 \, \dd \gl \;=\; \sum_{\ga \in \Sigma} d_\ga \,
\|\somb Q(\alpha) \|_{\text{hs}}^{2} \;<\; \infty \,.
\end{equation}
In particular, for every $\gep > 0$, 
$\|\somb Q(\alpha) \|_{\text{hs}} \le \gep$ for every $\ga \notin  \Gamma$,
with $\Gamma$ a finite subset of $\Sigma$.
 Since 
$\|\somb Q(\alpha) \|_{\text{op}} \le \|\somb Q(\alpha) \|_{\text{hs}}$, 
we have that $\|\somb Q(\alpha) \|_{\text{op}} \le \gep$ for every
$\ga \in \Sigma$, $\ga \not \in \Gamma$.
Next observe that Lemma~\ref{lem:1} can be applied, because
by hypothesis the support of $Q$ contains a non-empty open set $A$, hence
$A^{-1}A$ is open too and therefore it generates the whole $G$
(it is easily seen that, for any non-empty open subset $B$,
$\union_{n=1}^\infty B^n$ is non-empty and both open and closed, 
hence it must be the whole $G$, which is connected).
This observation yields
\begin{equation}
\label{eq:h}
	h \;:=\; \sup_{\ga \in \Sigma, \, \ga \ne 0} \|\somb Q(\alpha) \|_{\text{op}}
	\;<\; 1\,.
\end{equation}
Therefore from Lemma~\ref{lem:crucial} we have that
\begin{equation}
	\sup_\go \, E^\go(k) \;\le\; \|f\|_2 \cdot h^{(k-1)}\,,
\end{equation}
which proves \eqref{eq:toprove} under hypothesis H-\ref{H:2}.
\qed

\smallskip

\noindent
{\it Proof of \eqref{eq:toprove} under
Hypothesis H-\ref{H:1}.} 
Since the law $Q$ is assumed to be conjugation invariant, we have
$\int_G f(g) \, Q(\dd g) = \int_G f(t^{-1} g t) \, Q(\dd g)$, 
for every $t \in G$. Then for any law
$\nu$ on $G$ and for any bounded measurable function $f:G \to \R$ we have
\begin{equation}
\int_G f \, \dd (Q*\nu) \;=\; \int_G \int_G f(gh) \, Q(dg) \, \nu(dh) \;=\; \int_G \int_G
f(hg) \, Q(dg) \, \nu(dh) \;=\; \int_G f \, \dd (\nu*Q) \,,
\end{equation}
hence $Q * \nu = \nu * Q$. In particular, taking $\nu = \gd_g$,
the operator $\somb Q(\alpha)$ commutes with $U^\alpha(g)$, for every $g\in G$.
Schur lemma then yields that $\somb Q(\alpha)$ is a
multiple of the identity $I_\alpha$ on $H_\alpha$: $\somb Q(\alpha) =
c_\alpha I_\alpha$ for $c_\ga \in \C$.
Then from \eqref{eq:mydec} it follows that
\begin{equation}
\somb L^\go_{n+1,n+k}(\alpha)  \;=\;
U^\alpha(\go_{n+1} \cdots \go_{n+k}) \,
\somb Q(\alpha)^{k} \,,
\end{equation}
hence $\|\somb L^\go_{n+1,n+k}(\alpha)\|^2_{\mathrm{hs}}
= \|\somb Q(\alpha)^{k} \|^{2}_{\mathrm{hs}}$.
Since by assumption for $k \ge n_0$ the measure $Q^{*k}$ has a density
$f_k := \dd Q^{*k} / \dd \gl \in L^2(G,\dd \gl)$, 
it follows that also $L^\go_{n+1,n+k}(\alpha)$
has a density $g^\go_{n,k} = \dd L^\go_{n+1,n+k} / \dd \gl \in L^2(G,\dd \gl)$
(cf. Lemma~\ref{lem:Fourinv})
and by Parseval's identity we have
\begin{equation*}
	\int_G (g^\go_{n,k})^2 \, \dd \gl \;=\;
	\| g^\go_{n,k} \|_2^2 \;=\; \| f_k \|_2^2
	\;=\; \sum_{\ga \in \Sigma}
	d_\ga \, \big\| \somb{Q}(\ga)^k \big\|_{\text{hs}}^2 \;<\; \infty \,.
\end{equation*}
Arguing as above and recalling that $\somb Q(\ga) = c_\ga I_\ga$,
it follows that \eqref{eq:h} still holds.
We therefore have for $k \ge n_0$
\begin{equation}
\begin{split}
 & 4\, \|L^\go_{n+1,n+k}(\alpha) - \lambda\|^2_{TV} 
\;\le\; \left( \int_G |g^\go_{n,k} - 1| \dd \gl \right)^2
\;\le\; \int_G (g^\go_{n,k})^2 \dd \gl \,-\, 1 \\
& \;=\; \sum_{\alpha \in \Sigma, \, \ga \ne 0}
d_{\alpha} \, \big\|\somb Q(\alpha)^{k} \big\|^{2}_{\mathrm{hs}}
\;\le\; \sum_{\alpha \in \Sigma, \, \ga \ne 0}
d_{\alpha} \, \big\| \somb Q(\alpha)^{n_0} \big\|^{2}_{\mathrm{hs}} \,
\big\| \somb Q(\alpha) \big\|^{2 (k-n_0)}_{\mathrm{op}}  \;=\; \| f_{n_0} \|_2^2 \,  \cdot 
h^{2 (k-n_0)} \,.
\end{split}
\end{equation}
Then $\sup_\go E^\go(k) \le \| f_{n_0} \|_2 \, 
h^{k-n_0}$ and the proof
of equation \eqref{eq:toprove} is complete.
\qed

\bigskip

\noindent
{\it Proof of Proposition~\ref{th:mixing}.}
It suffices to prove that for every $n$ and every $\go$ we have
$\ga^\go (n) \, \le \,  2E^\go (n) $.
Since $\{\phi_{1,n}\}_{n}$ is a (inhomogeneous) Markov process 
we directly see that
\begin{equation}
    \ga^\go (n) \,\le \, 
    \sup_{u,w}  \left|\bE [u(\phi^\go_{1,m})w(\phi^\go_{1,m+n})]
    -\bE [u(\phi^\go_{1,m})]\bE[w(\phi^\go_{1,m+n})]
    \right|,
\end{equation}
where  $u$ and $w$ vary in the set of  measurable maps from
$G$ to $[0,1]$. Since
\begin{multline}
\left|\bE [u(\phi^\go_{1,m})w(\phi^\go_{1,m+n})]
    -\bE [u(\phi^\go_{1,m})]\bE[w(\phi^\go_{1,m+n})]
    \right|\, \le \\
    \left\vert
 \int_{G} u(g)\left(\int_G w(gg^\prime) \left(L_{m+1,m+n}^\go (\dd g^\prime)
   -\gl (\dd g^\prime)\right)\right) L_{1,m}^\go (\dd g) 
    \right\vert \,+\\
     \left\vert
     \int u(g) L_{1, m}^\go (\dd g)
 \int_{G} \left(\int_G w(gg^\prime) \left(L_{m+1,m+n}^\go (\dd g^\prime)
   -\gl (\dd g^\prime)\right)\right) L_{1,m}^\go (\dd g) 
    \right\vert
   ,
\end{multline}
the desired bound follows since both $\vert u(\cdot)\vert$ and $\vert w(\cdot)\vert$
are bounded by $1$.
\qed


\bigskip

\appendix

\medskip

\section{The elementary approach to the two-dimensional case}

\label{sec:2d}
We give here a partial proof of Theorem~\ref{th:main} in the 2-dimensional case. We identify in particular 
the variance $\gs^2$, {\sl cf.} \eqref{eq:sigma-def}, of the limit process. 
We set $\bbT := \R/(2\pi\Z)$ and we denote by $R_\ga$
the rotation by an angle $\ga$. With reference to \eqref{eq:X}, we write
$\go_j=R_{\gamma_j}$ and $r_j=R_{\theta_j}$, with $\gamma_j$ and
$\theta_j$ random variables taking values in $\bbT$. 
The Fourier coefficients of the law $Q$ of $\theta_1$ are 
$\int_\bbT e^{i m x} Q (\dd x)\, =: \hat{q}_m$, 
for $m\in\Z$. Recall that we are assuming that $Q(\{\theta_0, \theta_0+\pi\})<1$ 
for every $\theta_0$ 
and this is equivalent to $\vert \hat{q}_n\vert <1$ for $n=1$ and $n=2$.

 We set $\Theta_n := \theta_1 + \ldots + \theta_n$ and  $\Gamma_n := \gamma_1 +
\ldots + \gamma_n$ for $n\in\N$, along with $\Phi_n:= \Gamma_n +\Theta_n$.
Therefore the real and complex part of 
the random variable $Z_N^\go:= e^{i\Phi_1}+\ldots +
e^{i\Phi_N}$ coincide with the components
of the random vector $X_N^{v_0,\go}$, $v_0=(1,0)$.
Our goal is to compute the asymptotic covariance matrix of $X_N^{v_0,\go}$
as $N \to \infty$. Note that no centering is needed, since
\begin{equation}
    \bE\left[Z_N^\go\right] \;=\; \sum_{m=1}^N e^{i \Gamma_m} \, \bE\left[e^{i \Theta_m}\right]
    \;=\; \sum_{m=1}^N e^{i \Gamma_m} \, \hat{q}_1^m\,,
\end{equation}
and therefore $|\bE[Z_N]| \le |\hat{q}_1|/(1-|\hat{q}_1|) < \infty$,
because $|\hat q_1| < 1$.

We can therefore focus on the second moments.
For simplicity, we fix an arbitrary direction $e^{i\xi_0}$
in $\R^2 \simeq \C$, with $\xi_0 \in \bbT$,
and we look at the projection 
$\{Z^{\go, \xi_0}_n\}_n$ of the process $\{Z^\go_n\}_n$
in this direction, {\sl i.e.}
\begin{equation}
    Z^{\go, \xi_0}_0 \;:=\; 0\,  , \qquad \quad Z^{\go, \xi_0}_n \;:=\; \cos(\Phi_1 - \xi_0)
    \;+\; \ldots \;+\; \cos(\Phi_n-\xi_0)\,.
\end{equation}
For $n\in\N$ and $m\in \N\cup\{0\}$
one directly computes with $x=\Theta_m + \Gamma_{m+n}-\xi_0$
\begin{equation}
\label{eq:gsh1}
\bE \big[ \cos(\Phi_{n+m}-\xi_0) \,\big|\, \Theta_m \big] \,=\,
    \bE \big[ \cos(x+\Theta_n) \big] 
    \,=\, \Re \left( \hat q_1^n e^{ix}\right) \, =\, 
    \vert \hat q_1\vert ^n \cos(\Theta_m + \bar \theta \,n+ \Gamma_{m+n}-\xi_0)\,,
\end{equation}
where $\bar \theta$ is such that
$e^{i \bar \theta} = \hat q_1/\vert \hat q_1\vert$.
We observe also that for $a,b\in \bbT$
\begin{equation}
\label{eq:gsh2}
\bE \left[  \cos(\Theta_{m}+a)\cos(\Theta_{m}+b) \right] \,=\,
\frac 12 \cos (a-b) + \frac 12 \Re \left( \hat q_2^me^{i(a+b)}\right),
\end{equation}
and from \eqref{eq:gsh1} and  \eqref{eq:gsh2}  we directly see that
\begin{multline}
\label{eq:fromgsh}
\bE \left[  \cos (\Phi_m -\xi_0) \cos (\Phi_{n+m} -\xi_0)\right]\, =\\
\frac 12 \vert \hat q_1\vert^n 
\left\{ \cos (\Gamma _{m+n}-\Gamma_m +\bar \theta n) 
+ \Re \left( \hat q _2^m e^{i(\Gamma _{m+n} +
\Gamma_m +\bar \theta n -2\xi_0)}\right) \right\},
\end{multline}
and the latter expression actually holds also for $n=0$.
We are now ready to estimate $\bE[(Z_N^{\xi_0, \go})^2]$. 
The expression contains diagonal terms and for those we have 
\begin{equation}
    \sum_{m=1}^N \bE\big[ \cos^2(\Phi_m-\xi_0) \big] \;=\; \frac{N}{2}\; + o(N)\,,
\end{equation}
 by \eqref{eq:fromgsh} with $n=0$
(recall that $\vert \hat q_2 \vert <1$). 
The off-diagonal terms instead give
\begin{equation}
    2\sum_{m=1}^{N-1} \sum_{n=1}^{N-m} \bE\big[ \cos(\Phi_m-\xi_0)
    \cos(\Phi_{n+m}-\xi_0) \big] \;=\;
    \sum_{n=1}^{N-1} \hat q_1^n \sum_{m=1}^{N-n} \cos(\Gamma_{m+n} - \Gamma_m+\bar \theta n) \;+\; o(N)\,.
\end{equation}
For every fixed $n\in\N$, by the Ergodic Theorem we have that $\bbP(\dd \go)$--a.s. as $N\to\infty$
\begin{equation}
    \sum_{m=1}^{N-n} \cos(\Gamma_{m+n} - \Gamma_m+\bar \theta n) \;=\; \bbE(\cos(\Gamma_n+\bar \theta n))
    \cdot N \;+\; o(N) \,,
\end{equation}
and therefore that $\bbP$--a.s.
\begin{equation}
    \sum_{n=1}^{N-1} \hat q_1^n \sum_{m=1}^{N-n} \cos(\Gamma_{m+n} 
    - \Gamma_m + \bar \theta n) \;=\;
    \Bigg( \sum_{n=1}^\infty \hat q_1^n \,\bbE\left(\cos(\Gamma_n+\bar \theta n)\right) \Bigg) \cdot N \;+\; o(N)\,,
\end{equation}
so that finally we have $\bbP(\dd \go)$--a.s.
\begin{equation}
\label{eq:diffusion}
    \frac 1 N \bE\big[ (Z^{\go,\xi_0}_N)^2 \big]  \;=\;
    \frac{1}{N}\sum_{i,j = 1}^N \bE \big( \cos(\Phi_i-\xi_0)
    \cos(\Phi_j-\xi_0) \big) \stackrel{ N\to\infty }\longrightarrow
    \frac 12 \;+\; \sum_{n=1}^\infty \vert \hat q_1\vert^n \,\bbE[\cos(\Gamma_n +\bar \theta n)]\,,
\end{equation}
which matches with \eqref{eq:sigma-def}.
Note that the diffusion coefficient is independent of the direction
$\xi_0$ and that it depends on the law of
$\theta_1$ just through the first Fourier coefficient $\hat q_1$.


\medskip

\section{The homogeneous case}
\label{app:hom}

The aim of this appendix is to argue that, if disorder
is absent, Theorem~\ref{th:main} holds under the assumption
that the support of $Q$ generates a dense set in $SO(d)$. 

In order to do this, let us first observe that, 
when disorder is absent, we can weaken Assumption~\ref{ass:noise}
to the following {\sl generalized condition}: there exist $m \in \N$ such that
\begin{equation} \label{eq:homcond}
	\big\| \big( \bE(r_1) \big)^m \big\|_{\text{op}} \;<\; 1 \,, \qquad 
	\text{ and }\qquad
	\big\| \big( \bE(r_1 \otimes r_1) \big)^m \big\|_{H^0_s, \text{op}} \;<\; 1 \,,
\end{equation}
where we recall that $H^0_s$ denotes the space of symmetric real matrices
with zero trace.
We have shown in~\S\ref{sec:validity} that 
this condition with $m=1$ follows from Assumption~\ref{ass:noise}.
The fact that, when disorder is absent, equation \eqref{eq:homcond} is sufficient
to yield Hypothesis~\ref{hyp:main}, and hence Theorem~\ref{th:main}, is immediately
checked: for instance, by \eqref{eq:val1} we can write 
\begin{equation}
	\sum_{k=1}^\infty \big\| \bE[\phi^\go_{n,n+k}] \big\|_{\text{op}}
	\;\le\; \sum_{k=1}^\infty \big\| (\bE( r_1 ))^m \big\|_{\text{op}}^{\lfloor k/m \rfloor}
	\;<\; \infty\,,
\end{equation}
and analogously one shows that 
$\sum_{k=1}^\infty \big\| 
\bE\big[\,\overline{\phi^\go_{n,n+k} \otimes \phi^\go_{n,n+k}}\,\big]
- \Pi \big\|_{\text{op}} < \infty$, cf.~\eqref{eq:analogy3}.

We recall that, for an arbitrary linear operator $A$
on some vector space and for any fixed operator norm $\| \cdot \|$, 
the sequence $\| A^m \|^{1/m}$ converges as $m \to \infty$
toward the spectral radius of $A$, denoted $\text{Sp}(A)$. 
Furthermore, by sub-additivity 
(since $\|A^{m+n}\|_{\text{op}} \le \|A^m\|_{\text{op}} \, \|A^n\|_{\text{op}}$) 
we have
$\text{Sp}(A) =
\inf_{m \in \N} \| A^m \|^{1/m}$, hence we can restate \eqref{eq:homcond} as
\begin{equation} \label{eq:homcond2}
	\text{Sp}\big( \bE(r_1) \big) \;<\; 1 \,, \qquad \qquad
	\text{Sp}\big( \bE(r_1 \otimes r_1) |_{H^0_s} \big) \;<\; 1 \,.
\end{equation}

\smallskip

Let us finally show that equation \eqref{eq:homcond2} is satisfied whenever
the support $V$ of $Q$ generates a dense set in $SO(d)$, i.e., whenever
the closure of $\union_{k \in \Z} V^k$ is the whole $SO(d)$,
where we set $V^{-1} := \{g^{-1}: \, g \in V\}$, $V^2 := \{g \, h: \, g,h \in V\}$,
and so on. Since this fact is easily checked for $d=2$,
in the following we assume that $d \ge 3$.

We argue by contradiction: if the spectral radius of $\bE(r_1)$
is equal to one, there exists $v \in \C^d$ with $\|v\| = 1$ such that
$\bE(r_1)\, v = e^{i \theta}\, v$, with $\theta \in [0,2\pi)$, hence
\begin{equation}
	1 \;=\; \Re \langle e^{i\theta} v , \bE(r_1)\, v \rangle \;=\;
	\int_{SO(d)} \Re \langle e^{i\theta} v , g\, v \rangle \, Q(\dd g) \,.
\end{equation}
In the preceding relations we have denoted by $\langle \cdot, \cdot \rangle$
the standard Hermitian product on $\C^d$, i.e.,
$\langle a, b \rangle := \sum_{k=1}^d \overline{a_k} \, b_k$,
where $\overline a$ denotes the complex conjugate of~$a$.
Since $\Re \langle e^{i\theta} v , g\, v \rangle \le 1$ for every $g \in SO(d)$,
we must have $g \, v = e^{i\theta}\, v$ for every $g \in V$, the support of~$Q$.
Writing $v_1 + i\, v_2$ with $v_1, v_2 \in \R^d$ and denoting by $U$ the
linear subspace of $\R^d$ spanned by $v_1, v_2$, it follows that
$g \, U = U$ for every $g \in V$. Since by assumption $V$ generates
a dense set in $SO(d)$, by continuity we must have
$g \, U = U$ for every $g \in SO(d)$, which is clearly impossible
because $1 \le \dim(U) \le 2$ (recall that we assume $d \ge 3$).

With analogous arguments, if the spectral radius of
$\bE(r_1 \otimes r_1)$ on the space $H_s^0$ equals one, there must exist 
$v_1, v_2 \in H^0_s$ with $\| v_1 \|_{\text{hs}}^2 + \| v_2 \|_{\text{hs}}^2 = 1$
and $\theta \in [0,2\pi)$ such that
$g \, (v_1 + i\,v_2) \, g^{-1} = e^{i\theta} (v_1 + i\,v_2)$,
for every $g \in V$. Denoting by $U$ the
linear subspace of $H^0_s$ spanned by $v_1, v_2$, it follows that
$g \, U \, g^{-1} = U$ for every $g \in V$.
Since by assumption $V$ generates
a dense set in $SO(d)$, by continuity we must have
$g \, U \, g^{-1} = U$ for every $g \in SO(d)$.
However this is not possible, because 
the only linear subspaces $W$ such that $g \, W \, g^{-1} \subseteq W$
for every $g \in SO(d)$ are $W=\{0\}$ and $W = H^0_s$
(i.e., the representation $SO(d) \ni g \mapsto g \otimes g$ 
on the vector space $H^0_s$ is {\sl irreducible}).

Let us check this fact. We take $w \in W$ not
identically zero: by the spectral theorem, there exists $g \in SO(d)$ such
that $v := g \, w \, g^{-1} \in W$ is diagonal: $v_{ij} = \gl_i \, \gd_{ij}$.
Since $v$ is not identically zero and it has zero trace, there exist
$i_0, j_0$ such that $\gl_{i_0} \ne \gl_{j_0}$. Let us now take $h \in SO(d)$
to be the matrix that permutes the coordinates $i_0$ and $j_0$, i.e.,
$h_{ij} := \gd_{ij}$ for $i,j \not \in \{i_0, j_0\}$ 
while $h_{i_0 j} = h_{j i_0} := \gd_{j_0 j}$
and $h_{i j_0} = h_{j_0 i} := \gd_{i i_0}$. It is clear that
$\tilde v := h \, v \, h^{-1} \in W$ 
is such that $\tilde v_{ij} = \tilde \gl_i\, \gd_{ij}$,
where $\tilde \gl_i = \gl_i$ for 
$i \not \in \{i_0, j_0\}$ while $\tilde \gl_{i_0} = \gl_{j_0}$
and $\tilde \gl_{j_0} = \gl_{i_0}$. Therefore 
$z :=\frac{1}{(\gl_{i_0} - \gl_{j_0})} (v - \tilde v) \in W$ is such that
$z_{i_0 i_0}=1$, $z_{j_0 j_0}=-1$, and $z_{ij}=0$ for all the other values of $i,j$.
By considering $g \, x \, g^{-1}$, where $g \in SO(d)$ is an 
arbitrary permutation matrix,
we obtain all the matrices defined like $z$ but with arbitrary $i_0, j_0$.
These matrices span the linear subspace consisting of all the diagonal matrices
with zero trace, which are therefore contained in $W$. However, again by the spectral
theorem, for any matrix $u \in H^0_s$ we can find $g \in SO(d)$ such that
$g \, u \, g^{-1}$ is diagonal with zero trace, hence we must have $W = H^0_s$
and the proof is completed.

\medskip

\section*{Acknowledgments}
We acknowledge the support of ANR, grant POLINTBIO.


\end{document}